\documentclass[12pt]{article}

\usepackage{graphicx}
\usepackage{enumerate}
\usepackage{natbib}

\usepackage{graphicx,bm,url,microtype}
\usepackage{paralist}
\usepackage{authblk}
\usepackage{amsfonts,amssymb,amsthm,amsmath, amsthm}

\usepackage[a4paper,text={16.5cm,25.2cm},centering]{geometry}
\usepackage[compact,small]{titlesec}

\newtheorem{proposition}{{\sc\bf Proposition}}
\newtheorem{theorem}{{\sc\bf Theorem}}
\newtheorem{lemma}{{\sc\bf Lemma}}

\newtheorem{corollary}{{\sc\bf Corollary}}

\newtheorem{example}{{\sc\bf Example}}

\usepackage{amsmath,amsthm,mathrsfs,amsfonts}
\usepackage{ulem} 	
\usepackage{url}
\usepackage{amsfonts,todonotes}
\usepackage{amssymb}
\usepackage{bbm}
\usepackage{graphicx}
\usepackage{subfig}
\usepackage{float}

\def\defeq{\mathrel{\mathop:}=}

\newcommand{\I}{\mathbbm{1}}

\usepackage{multirow}
\newcommand*{\QEDB}{\hfill\ensuremath{\square}}%

\usepackage{parskip}
\usepackage{booktabs}
\usepackage[shortlabels]{enumitem}
\usepackage{unicode}

\begin{document}

\begin{center}
	\Large \bf Weighted lens depth: Some applications to supervised classification
\end{center}
\begin{center}
		Alejandro Cholaquidis $^1$, \hspace{.2cm}
		Ricardo Fraiman $^1$ \\ Fabrice Gamboa$^2$ \hspace{.2cm}  and \hspace{.2cm}  Leonardo Moreno$^3$ \\
		$^1$ Centro de Matemáticas, Facultad de Ciencias, Universidad de la República, Uruguay. \\
		$^2$ Institut de Mathématiques de Toulouse, France. \\
		$^3$ Instituto de Estadística, Departamento de Métodos Cuantitativos, FCEA, Universidad de la República, Uruguay.
\end{center}

\begin{abstract}
Starting with Tukey's pioneering work in the 1970's, the notion of depth in statistics 
has been widely extended especially in the last decade.
These extensions include high dimensional data, functional data, and manifold-valued data.
 In particular, in the learning paradigm, the depth-depth method has become a useful technique.
In this paper we extend the notion of lens depth to the case of data in metric spaces, and prove its main properties, with particular emphasis on the case of Riemannian manifolds, where we extend the concept of lens depth in such a way that it takes into account non-convex structures on the data distribution.
Next we illustrate our results with some simulation results and also in some interesting real datasets, including pattern recognition in phylogenetic trees using the depth--depth approach.
\end{abstract}

\section{Introduction}

The notion of depth was first introduced in statistics by Tukey in the 70's, see \cite{tukey1975}, for multivariate data.
It is defined on $\mathbb{R}^d$ by the Tukey depth function $\textrm{TD}(x,P): \mathbb R^{d} \to [0, \infty)$ given by
$
\textrm{TD}(x, P)\defeq \inf \{ P(H): \textrm{$H$ a closed halfspace}, x \in H \},$  for a probability measure $P$ on $\mathbb{R}^d$. Roughly speaking, given $P$, $\textrm{TD}(x,P)$ measures how deep $x$ is with respect to  $P$.
Following that seminal paper, several notions of depth were proposed and developed.
 Two other important and well established notions of depth are the simplicial depth and  the  spatial depth.
The simplicial depth, introduced by Liu, see \cite{liu1990, liu1992, liu1992b},
is defined by $
\textrm{SD}(x,P_X)\defeq  \mathbb{P} \left( x \in S[X_1, \ldots,X_{d+1}]  \right),$  where $X_1, \ldots,X_{d+1}$ is an iid sample of $\mathbb{R}^d$ with common distribution $P_X$, and $S[X_1, \ldots,X_{d+1}]$ is the $d$-dimensional simplex with vertices 
 $X_1, \ldots,X_{d+1}$, that is, the set of all possible convex combinations of the sample points. Spatial depth, introduced by Serfling, see \cite{serfling2002,vardi2000}, is defined by $\textrm{SpD}(x,P) = 1 -  \mathbb E_{P_X}(S(x-X)),$ where for $x\neq 0$, $S(x)$ denotes the projection on the unit sphere and $S(0)=0$.
Several other depth measures have been proposed for different kinds of data, for instance, \textit{convex hull peeling depth} \cite{barnett1976}, \textit{Oja depth} \cite{oja1983},  and \textit{spherical depth} \cite{elmore2006}, among others.
 Some of them can be easily generalized to infinite-dimensional settings, like Tukey's depth.
However, this is not the case of some other ones, like  Liu's.
In this setting,  other notions of depths have been developed, see for instance   \cite{fraiman2001,lopez2009,claeskens2014,cuevas2009}.
More recently, notions  of depth have also been studied in general, non-Euclidean, metric spaces, see \cite{fraiman2019}.
We refer to \cite{serfling2000},  which discuss properties that a depth measure should have (like, for instance, vanishing at infinity).
However, some of these properties can not be extended to general metric spaces.

In this paper we will focus on the \textit{lens depth} (denoted by LD), which was introduced in \cite{liu2011}, in $\mathbb{R}^d$,
\begin{equation}
\label{ld}
\textrm{LD}(x,P_X)=  \mathbb{P} \left( x \in B \left(X_1, \Vert X_2-X_1 \Vert \right) \cap  B \left(X_2, \Vert X_2-X_1 \Vert  \right)  \right), \quad x\in \mathbb{R}^d.
\end{equation}

\noindent Here,  $B(x,r)$ denotes the closed ball centred at $x$ with radius $r$, and,  as before, $X_1$, $X_2$ are  independent random variables with common distribution $P_X$.
It enjoys many good properties: for example, its computational complexity in $\mathbb{R}^d$ is of order $O(n^{2}d)$, with $n$ being the sample size (see \cite{elmore2006}).  This is considerably less than the empirical computational complexity of Tukey's depth, $O\left(n^{(d-1)}\right)\log(d)$, or Liu's one, $O\left(n^{(d+1)}\right)$.

Since this depth depends only on the distances between points (or any measure of similarity between them), it can be easily extended to any metric space $(M,\rho)$ by replacing the norm in (\ref{ld}) with the distance $\rho$.
In particular, this is the case when $M$ is a Riemannian manifold and $\rho$ is the geodesic distance.
   We will also prove, in this very general setting, that this depth has some of the desirable properties stated by Serfling.
The empirical version of this depth for an iid sample $X_1, \ldots, X_n$ with distribution $P_X$ is obtained by replacing $P_X$ by the empirical distribution of the data, $P_n$.
We show the consistency of the empirical version with the population one, and also provide the convergence rate. Another aim of our work is to introduce and study, for Riemannian manifolds, a more general version of this depth, which will be called the weighted lens depth.
We will denote it by $\textrm{WLD}_p$.
To build this depth measure, the idea is to take into account the underlying structure of the data, provided  by the density of $P_X$, by changing the Riemannian metric on the manifold. 
 For this generalized version we will also provide the convergence rate for the empirical version.
It is important to note that for Riemannian manifolds embedded in $\mathbb{R}^d$, this empirical version can be computed and is consistent even if the Riemannian structure is unknown.

The classification problem for data on manifolds has recently gained importance, see for instance  \cite{yao2020} and the references therein. With a different approach, we will apply  $\textrm{WLD}_p$ on manifolds, to a classification method called the depth-depth, see \cite{liu1999, li2012, cuesta2017}, and show the performance of the depth-depth method through the study of two real-life data sets.

Our construction enjoys three very nice properties. \textbf{Computational complexity}: Some depths (like Liu depth) are computationally infeasible in high dimensional problems, because their complexity grows exponentially with the dimension.
This is not the case of lens depth. \textbf{Generality}: Lens depth (and also $\textrm{WLD}_p$) easily extends to general metric spaces and in particular to Riemannian manifolds.	 \textbf{Level-set structure}: Some depths (under restrictive conditions like unimodality or symmetry) characterize the distribution of the data, \cite{kong2010,kotik2017}, and have convex level sets.
 However, if the aforementioned  very restrictive  conditions are not fulfilled, the level sets are not convex, see \cite{dutta2011}.
Except for the case of distributions with convex support (see \cite{hlubinka2013}), the classical notions of depth are not designed to capture the underlying structure of the distribution.
For instance, the balls defining the lens depth are balls with respect to the norm, and they do not take into account the density of the data.
This makes the aforementioned classification method perform badly.
To overcome this issue, several modifications of the classical notions of depths have been proposed.
For instance in \cite{kotik2017} a weighted Tukey's depth is considered.

As shown in \cite{li2012}, the depth-depth approach  behaves better than other learning procedures when the data have some non--standard patterns.

\begin{figure}[h] 
\centering
\includegraphics[width=85mm]{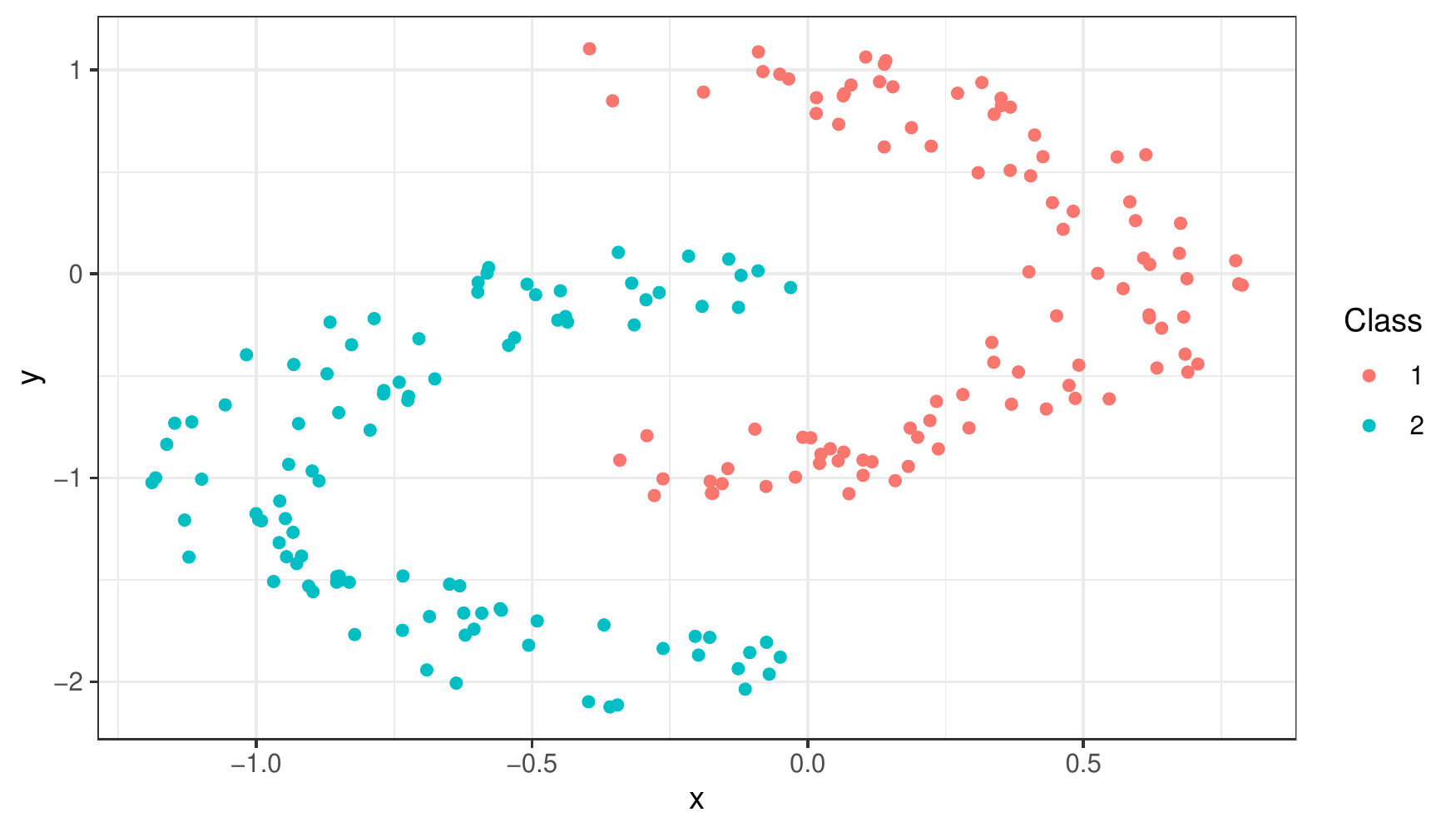}
\caption{200 data points  for Example 1.}
\label{medialunas}  
\end{figure}

\begin{figure}[!ht] 
\centering
\subfloat{\includegraphics[width=70mm]{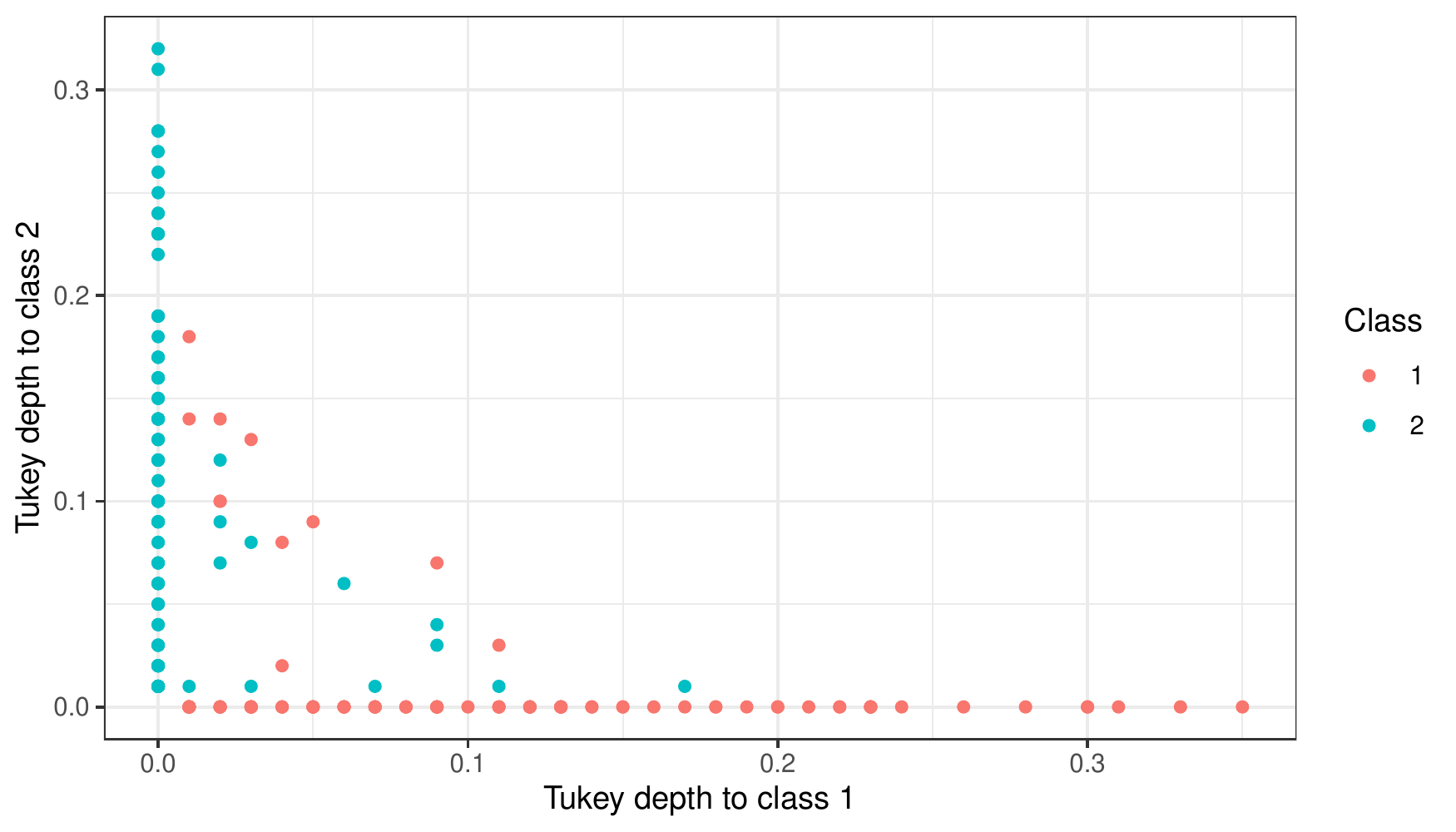}}
\subfloat{\includegraphics[width=70mm]{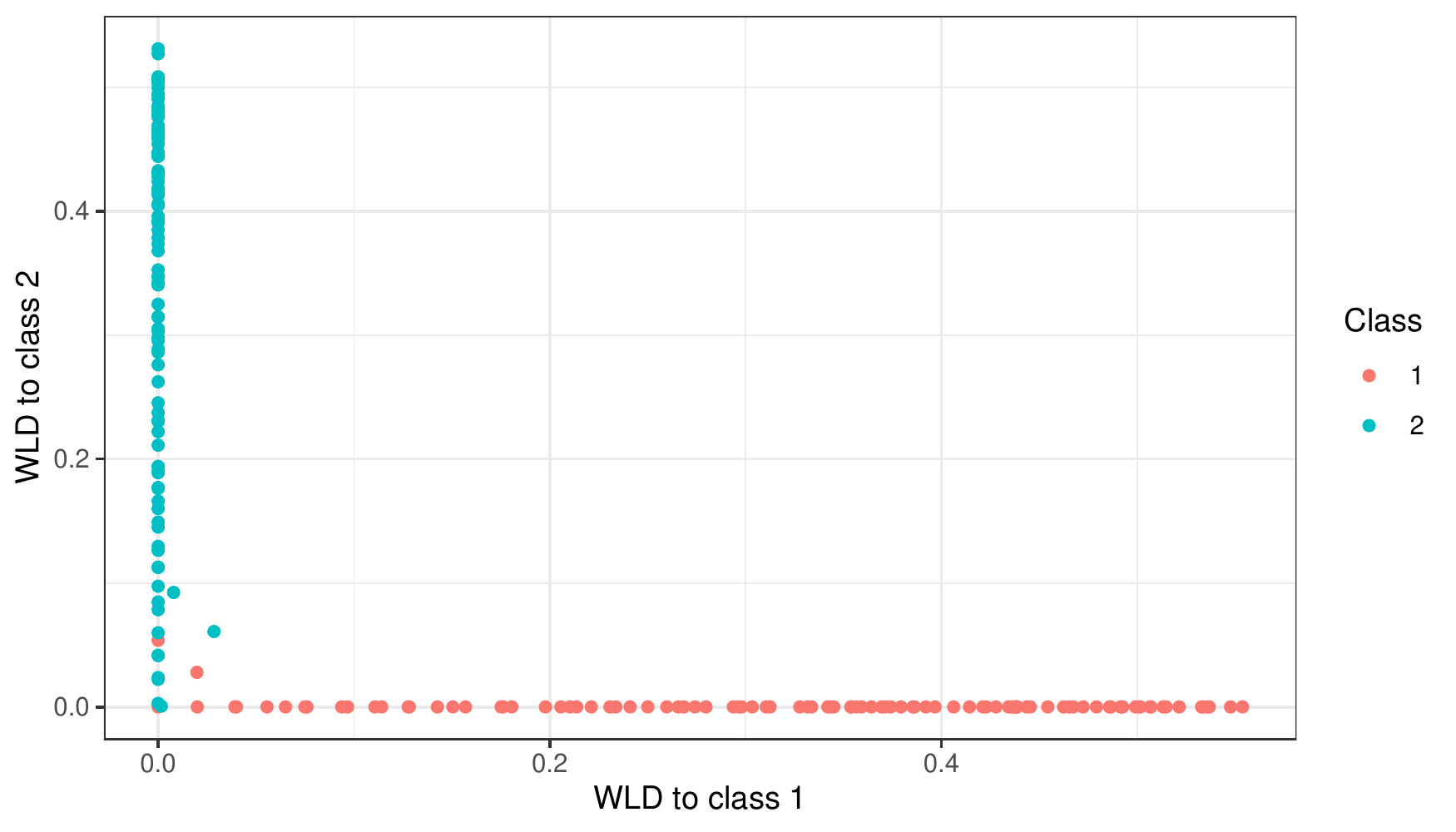}}
\caption{Depth-depth plots for the data shown in Figure \ref{medialunas}.
Left Panel: Results using Tukey's depth.
Right Panel: Results using $\textrm{WLD}_p$ with $p=2$.} \label{dep11}  
\end{figure}

The first example of the application of the aforementioned depth-depth technique (see  \cite{cuesta2017}) will be a classification problem where we have two classes, both of  sample size $100$.
The data points are shown in Figure \ref{medialunas}.
 We choose at random 50 of the 200 original sample points to use as a test sample.
 We apply Random Forest, trainned with the pairs $(\widehat{\textrm{TD}}_0(X_i), \widehat{\textrm{TD}}_1(X_j))_{i,j}$, for $i,j=1,\dots, 150$, where $\widehat{\textrm{TD}}_1(X_i)$ is Tukey's depth of $X_i$ w.r.t. the 1-class, and  $\widehat{\textrm{TD}}_2(X_i)$ is Tukey's depth of $X_i$ w.r.t. the 2-class, see Figure \ref{dep11}.
The misclassification error rate over 1000 replications was $10\%$.
If we do the same but using $\textrm{WLD}_p$ with $p=2$, the misclassification error rate obtained was $1\%$.

\subsection{Outline}

This paper is organized as follows.
Section \ref{ld1} extends the concept of lens depth to metric spaces and proves some of its properties.
Section \ref{wld1} introduces the concept of weighted lens depth and its estimator for  data valued in a Riemannian manifolds.
Consistency as well as convergence rates are obtained.
Section \ref{gg} uses $\textrm{WLD}_p$ in a  supervised classification by means of the depth-depth method,  for simulated data and  for examples of real-life data.
All proofs are postponed to the Appendix.

\section{Notation}

In what follows, $(M,d)$ is a complete separable metric space, endowed with the Borel $\sigma$-algebra.
The  closed ball  centred at $x$ of radius $r>0$ is denoted by $B(x,r)$.
 We denote by  $X$ (or more generally for $j\in \mathbb{N}^*$, $X_j$) an $M$-valued random variable, whose distribution is denoted by $P_X$.
Given a set $A\subset M$, $\partial A$ denotes its boundary.
Given $\epsilon>0$, we set $A^\epsilon=\{x\in M: d(x,A)\leq \epsilon\}$ (sometimes also denoted by $B(A,\epsilon)$).
For $x\in M$,  we write $d(x,A)=\inf_{a\in A} d(x,a)$.
Given two non-empty compact sets $A,C$, their Hausdorff distance $d_H(A,C)$  is defined by
\begin{equation}\label{hausdist}
	d_H(A,C)=\max\Big\{\max_{a\in A}d(a,C), \ \max_{c\in C}d(c,A)\Big\}.
\end{equation}
 Given $x,y\in M$, we  set $A(x,y)=B(x,d(x,y))\cap B(y,d(x,y))$.

\section{Lens depth on metric spaces}\label{ld1}

Recall that we have defined in (\ref{ld}), for $x\in M$, the lens depth of $x$ to be $\textrm{LD}(x,P_X)=\mathbb{P}(X\in A(X_1,X_2)).$ Given an iid sample $\mathcal{X}_n= \{X_1, \ldots, X_n\} \subset M$ of a random variable $X$ taking values in $M$, the empirical version of LD is defined by 
\begin{equation}
	\widehat{\textrm{LD}}_n(x) = \binom{n}{2}^{-1} \sum_{1 \leq i_1 < i_2 \leq n} \I_{A(X_{i_1} X_{i_2})}(x).
	\label{estimadorprofundidad_cape}
\end{equation}
Let us observe that  $\widehat{\textrm{LD}}_n(x)$ is a U-statistic of order $2$.
In the next section, we extend, to general metric spaces, some properties of LD and $\widehat{\textrm{LD}}_n(x)$ which hold in the finite dimensional setting.

\subsection{Some properties}

In this section we will recall some desirable properties of the depth LD (see \cite{serfling2000} or \cite{fraiman2019}), for metric spaces.
To prove the continuity of LD, we will make an assumption on the random variable $X$.
Namely, we assume 
\begin{equation}\label{lconti}
\mathbb{P} \big( d(X,x_1)= d(X,x_2) \big)=0, \quad  \textrm{for all  }  x_1,x_2\in M, x_1\neq x_2.
\end{equation}
In $n$-dimensional Euclidean space, this condition is equivalent to assuming that the  $(n-1)$-dimensional linear subspaces  have probability zero.
If condition \ref{lconti} holds, then given $x \in M$, 
\begin{align}\label{lcont}
 \mathbb{P}  \big ( d(X_1,x) = d(X_1,X_2)  \big) = 
  \int_M \mathbb{P} \big ( d(X_1,x) = d(X_1,t)  \big) dP_X(t)=0.
\end{align}

\begin{proposition} {(Vanishing at infinity).}
\label{prop1}
Let $(M,d)$ be a separable and complete metric space such that there exists $z_0$ with $\sup_{t\in M} d(z_0,t)=\infty$.
Let $y\in M$ be fixed.  Then,

\begin{enumerate}
	\item[(a)] $\sup_{d(x,y)>K} {\textrm{LD}(x)} \rightarrow 0 $   as  $K \rightarrow +\infty,$  and
	\item[(b)] $ \lim_{K \to \infty} \sup_{d(x,y)>K} \widehat{{LD}}_n(x) =0$ a.s.
\end{enumerate}
\end{proposition}

\begin{proposition} {(Continuity and point-wise consistency).}
\label{prop2}
Let $x,y \in M$ and assume that  $X$ satisfies \eqref{lconti}.
Then, 
\begin{enumerate}
	\item [(a)](Continuity) $\vert {\textrm{LD}}(x) -  {\textrm{LD}}(y) \vert \rightarrow 0  \quad \textrm{when} \quad d(x,y) \rightarrow 0$.
	\item [(b)](Consistency-rate of convergence) $\beta_n|\widehat{\textrm{LD}}_n(x)-{\textrm{LD}}(x)|\to 0$ a.s.
for all $\beta_n\to\infty$ such that $n/(\beta_n^2\log(n))\to\infty$.
\end{enumerate}
\end{proposition} 

\section{Weighted lens depth on manifolds}
\label{wld1}
\subsection{Basic concepts}
 
We will start by introducing briefly some basic concepts of Riemannian geometry,  that will be used throughout this section. A Riemannian manifold is defined by a differentiable manifold $M$ of dimension $d$, equipped with a Riemannian metric $g_1$ which defines for every point $p \in M$  the scalar product of tangent vectors in the tangent space $T_p M$, smoothly depending on the point $p$. 
 
Given two points $x,y \in M$  the \textit{induced distance} between them is defined as the infimum of the lengths $L(\gamma)$ of all  continuously differentiable curves   $\{\gamma_t\}_{t\in [0,1]}\subset M$ joining $x$ and $y$ .

  A geodesic (with speed $s \in \mathbb{R}^{+}_0$) is a smooth map $\alpha:I \rightarrow M$, where $I$ is an interval, such that $\Vert \alpha^{'}(t) \Vert= s$  for all $t \in I$ and which is ``locally length minimizing''.
 

Through  this section we  assume that $(M,g_1)$ is a $d$-dimensional Riemannian manifold without boundary (which is assumed to be unknown), and the available data is an iid sample $\mathcal{X}_n= \{X_1,\dots,X_n\} \subset M$ of a random vector  $X$ whose distribution $P_X$ is assumed to be supported on $M$. We assume it has a density $f$ w.r.t. the volume form on $M$ given by the metric $g_1$, see \cite{petersen2006}.
In this setting, to compute the empirical lens depth (defined by \eqref{estimadorprofundidad_cape}) requires computing the geodesic distance between two points.
This is  a problem that has been previously tackled in the literature (see for instance  \cite{cormen2009} and  \cite{davis2019}).
We will introduce a more general version of the lens depth, which will be called the weighted lens depth, that takes into account the density $f$ together with the underlying Riemannian structure.

Given a parameter $p\geq 1$ (called the power parameter), we define a new metric $g_p= f^{2(p-1)/d}g_1$, see \cite{hwang2016}.  The geodesic distance w.r.t. $g_p$ between $x,y\in M$, is denoted by $\rho_p(x,y)$. Let $A\subset M$ be a locally finite set, that is  $A\cap B$ is a finite set for any $B\subset M$ of finite volume w.r.t. $g_1$. We define $L_p(x,y;A)$ by  
$$L_p(x,y;A)= \min_{ \{ x_0=x,\ldots,x_k=y \}  \subset A^* } \sum_{i=0}^{k-1} (\rho_1 (x_i, x_{i+1}))^p, $$
where $A^*= A \cup \{x,y\}$.
We denote by $B_{p}(x,r)$ the closed ball of radius $r$, w.r.t. $\rho_p$, centred at $x$.
 We have
\begin{equation}\label{fermat}
\rho_p(x,y)= \inf_{\gamma} \int_0^1 f(\gamma_t)^{(1-p)/d} \sqrt{g_1(\gamma^{'}_t, \gamma^{'}_t)}dt,
\end{equation}

 where  $\{\gamma_t\}_{t\in [0,1]}\subset M$ is a piece-wise $C^\infty$-curve such that $\gamma_0=x$ and $\gamma_1=y$. The distance define by \eqref{fermat} is called Fermat distance. As mentioned in \cite{groisman2018}, ``it coincides with Fermat principle in optics for the path followed by light in a non-homogeneous media when the refractive index is given by $f^{-\beta}$".

In \cite{hwang2016} the following result is proven.  It gives the (asymptotic) relation between $\rho_p$ and $L_p(x,y;\mathcal{X}_n)$.

\begin{theorem}{(Hwang--Damelin--Hero, 2016)}
\label{teo1}
 Assume that the density $f$ of $X$ is continuous and bounded from below by a positive constant.
Then, there exists a  constant $C(d,p)>0$ such that for all $\epsilon>0$  and $b>0$, there exists a $\theta_0>0$ such that, for  $n$ large enough.
$$\mathbb{P} \Big( \sup_{x,y \in M, \rho_1(x,y) \geq b}
\left \vert \frac{L_p(x,y; \mathcal{X}_n)}{ n^{(1-p)/d}\rho_p(x,y)} -C(d,p)  \right \vert > \epsilon \Big) \leq \exp \left(- \theta_0 n^{1/(d+2p)}\right),$$

\end{theorem}
 
If $p=1$, it can be proved that $C(d,p)=1$ (see \cite{hwang2016}).
It is important to note that neither $\theta_0$ nor $C(d,p)$ depend on $n$.
A similar version of the previous theorem is stated in \cite{mckenzie2019}, 
with weaker assumptions, (the geodesic distance $\rho_1$ is replaced by the Euclidean norm $\Vert x_{i+1} - x_i \Vert$).
They also adress the clustering problem on manifolds.
If the manifold is unknown, then $\rho_1$ in the definition of $L_p(x,y;A)$ can be replaced by the Euclidean distance, as in \cite{groisman2018} or in \cite{mckenzie2019}. We extend our previous notations by defining the set  $A_p(x_1, x_2)= B_p(x_1, \rho_p(x_1,x_2)) \cap B_p(x_2, \rho_p(x_1,x_2))$ 
whose empirical version is $A_{p,n}(x_1, x_2)= A^1_{p,n}(x_1, x_2) \cap  A^2_{p,n}(x_1, x_2),$   where $ A^i_{p,n}(x_1, x_2) =  \left \{ t \in M :   L_p(t,x_i, \mathcal{X}_n) <  L_p(x_1,x_2, \mathcal{X}_n) \right \} \quad i=1,2.$
As will be shown in the following example, these sets depend heavily on the choice of $p$.
\begin{example} To see how the shape of the level sets of $WLD_p$ varies with $p$, let us consider the density 
	$f(u,Ku)=  C \exp \{ - \left((u-u_0)\alpha^{-1}\right)^2 -  \left( (Ku-Ku_0)\beta^{-1} \right)^2  \}, \quad  \textrm{with } u, Ku \geq 0,$
 where $u_0=0.2$, $\beta=0.5$,  $Ku_0=30$, $\alpha=3$ and $C$ is the normalization constant.
In Figure \ref{bolass}, the densities of the orthogonal and parallel components are shown, i.e.  $x=u Ku$ and $y=\sqrt{u^2-x^2}$, based on a sample of size 200.
This density has important applications when modelling  \textit{electron-cyclotron-maser (ECM) instability and models for electron energization}, see \cite{tong2017}.
As can be seen in  Figure \ref{bolass}, the sets $A^i_{2,200}(x_1, x_2)$  adapt their geometry to the level sets of the density.
This effect does not appear in  $A^i_{1,200}(x_1, x_2)$, since for $p=1$ this set depends only on the extremes $x_1,x_2$.
\begin{figure}[!ht] 
\centering
\subfloat{\includegraphics[width=55mm]{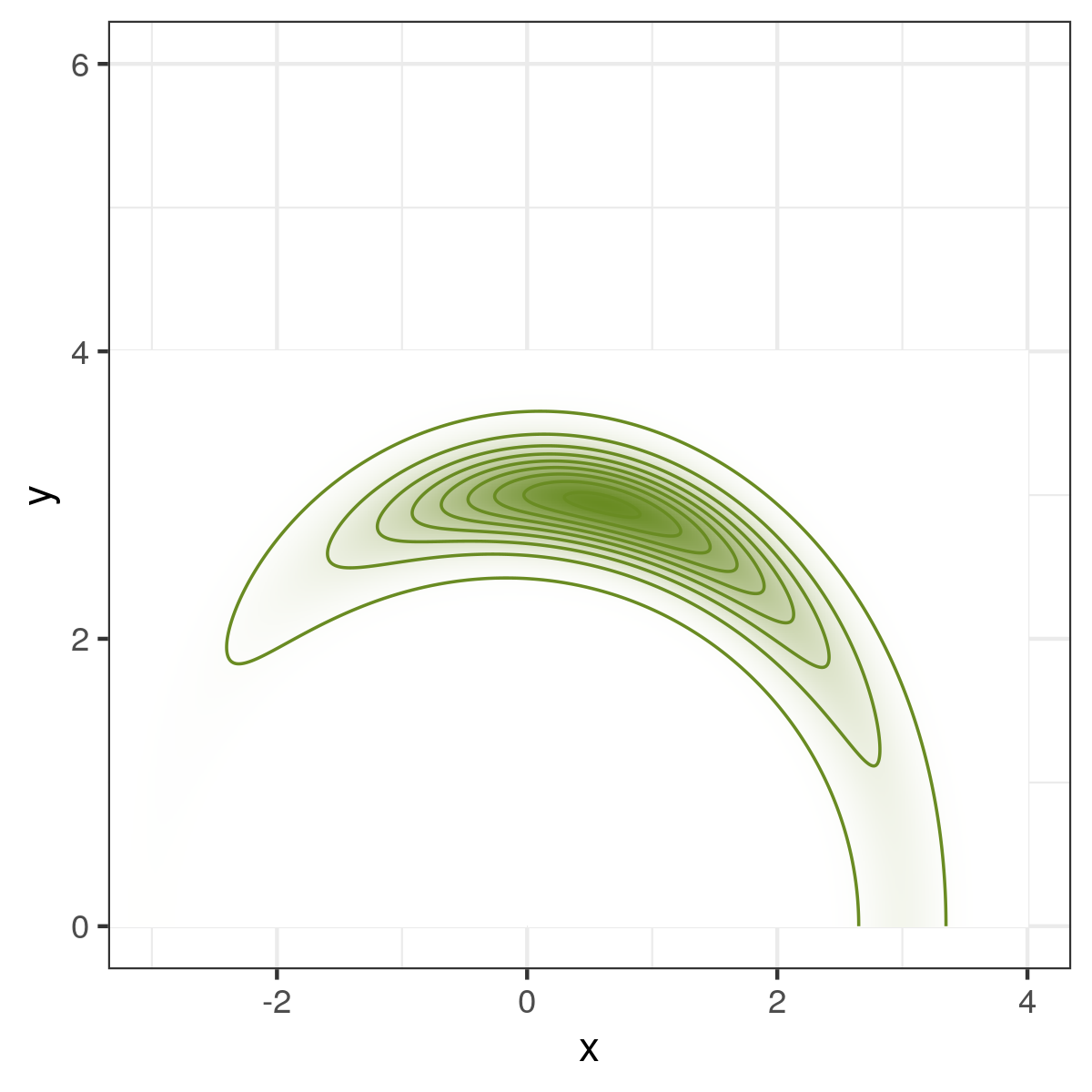}}
\subfloat{\includegraphics[width=55mm]{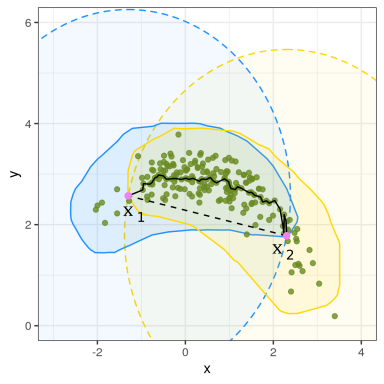}}
\caption{Left panel: Level sets of the density.
Right panel: The sets $A^1_{1,n}(x_1, x_2)$ and $A^2_{1,n}(x_1,x_2)$ in blue, while  $A^1_{2,n}(x_1, x_2)$ and $A^2_{2,n}(x_1,x_2)$ are shown in yellow.
Sample points are shown in green.
The black dotted line shows the geodesic joining $x_1$ and $x_2$ for $p=1$. The piece-wise black curve is the geodesic joining $x_1$ and $x_2$ for $p=2$.} \label{bolass}  
\end{figure}
\end{example}
\subsection{Main results}

 With the previous notations we define the weighted lens depth of a point $x\in M$, of order $p$, denoted by $\textrm{WLD}_p(x)$, by  
\begin{equation*}
\begin{aligned}
\textrm{WLD}_p(x) & = \mathbb{P} \left( x \in A_p(X_1,X_2) \right)
& = \int_{M^2} \I_{A_p(x_1,x_2)}(x) P_X(dx_1) P_X(dx_2).
\end{aligned}
\end{equation*}

The empirical version of $\textrm{WLD}_p$, is given by the  following  $U$-statistic of order two: 
\begin{equation*}
\widehat{\textrm{WLD}}_{p,n}(x) \defeq \binom{n}{2}^{-1} \sum_{1 \leq i_1 < i_2 \leq n} \I_{A_{p,n}(X_{i_1} X_{i_2})}(x).
\end{equation*}

\subsection{Consistency of the estimator}

In order to prove that $\widehat{\textrm{WLD}}_{p,n}(x)$ is a consistent estimator, we first prove a lemma about the Hausdorff distance between the sets $B_p(x_1,\rho_p(x_1,x_2))$ and its empirical counterpart $A^1_{p,n}(x_1,x_2)$, where the Hausdorff distance  $d_{H}$ is computed using \eqref{hausdist} for $d=\rho_1$. 
\begin{lemma} 
\label{lema1}
 Under the assumptions of Theorem \ref{teo1},
let $x_1,x_2 \in M$  and $\delta>0$, there exists $\theta_0>0$ such that, for $i=1,2$, $\mathbb{P}  \big( d_{H} \left(B_p(x_i, \rho_p(x_1,x_2)), A^i_{p,n}(x_1, x_2) \right) > \delta \big ) \leq  \exp \left(- \theta_0 n^{1/(d+2p)}\right), $ for $n$ large enough,  where $\theta_0$ is the constant given in Theorem \ref{teo1}.
\end{lemma} 

As a direct consequence of this lemma, we have the following corollary.
\begin{corollary}
\label{c1} 
\begin{itemize}
\item[]
\item[1)] For $i=1,2$, $d_{H} \left( B_p(x_i, \rho_p(x_1,x_2)), A^i_{p,n}(x_1, x_2) \right)\to 0$ a.s., as  $n \rightarrow \infty$.
\item[2)] For $n$ large enough, $P(d_{H} \left( A_p(x_1, x_2) , A_{p,n}(x_1, x_2) \right)>\delta)\leq 2\exp \left(- \theta_0 n^{1/(d+2p)}\right).$
\end{itemize}
\end{corollary}

\textbf{Condition I.} We will say that a manifold $(M,g_1)$ fulfils condition I if for all $p\in M$, $\partial B_1(p,\gamma)$ has positive reach for almost all $\gamma\in\mathbb{R}^+$ w.r.t. the Lebesgue measure on $\mathbb{R}$.

 Recall that the reach of a set $A \subset (M,d)$ is defined as
$
\sup_{r \in \mathbb{R}} \{ \forall x \in M, \ \mbox{for which} \ d(x, A)< r, \ \mbox{there exists a unique point} \ y \in M \ \mbox{such that} \  d(x,y)= d(x, A)\}.
$

 Theorem 5.3 in \cite{rataj} states that condition $I$ is fulfilled in dimensions 2 and 3 if the manifold $M$ is complete.
 
\begin{theorem}\label{conswld} Under the hypotheses of Theorem \ref{teo1}. Assume that $M$ fulfills condition I. Assume further that the density $f$ is $C^2$. 

Then, for all $x\in M$, 
\begin{equation}\label{rate}
\beta_n|\widehat{\textrm{WLD}}_{p,n}(x)- \textrm{WLD}_p(x)|\to 0\quad a.s., \quad \text{ as } n\to\infty, \ \mbox{ where }\ n/(\beta_n^2 \log(n))\to \infty.
\end{equation}

\end{theorem}

\section{Level sets of $\textrm{WLD}_p$ in $\mathbb{R}^2$}
\label{simfin}
We start by considering two distributions on the plane: a bivariate exponential distribution and an uniform distribution on a ring. 
We will study the shapes of the level sets when varying $p\in \{1,1.5,2,5\}$,  for a sample of size   $n=400$.
Let us recall that for $p=1$, we obtain the usual lens depth with the Euclidean distance.
Two main scenarios are considered. (a) Bivariate exponential, with density $f(x,y)= e^{-x-y}$  with  $x,y >0.$. (b) Uniform distribution on the ring $B((0,0),1.5)\setminus B((0,0),1)$.

\begin{figure}[ht]
 \centering
    \subfloat[]{\includegraphics[width=0.46\textwidth]{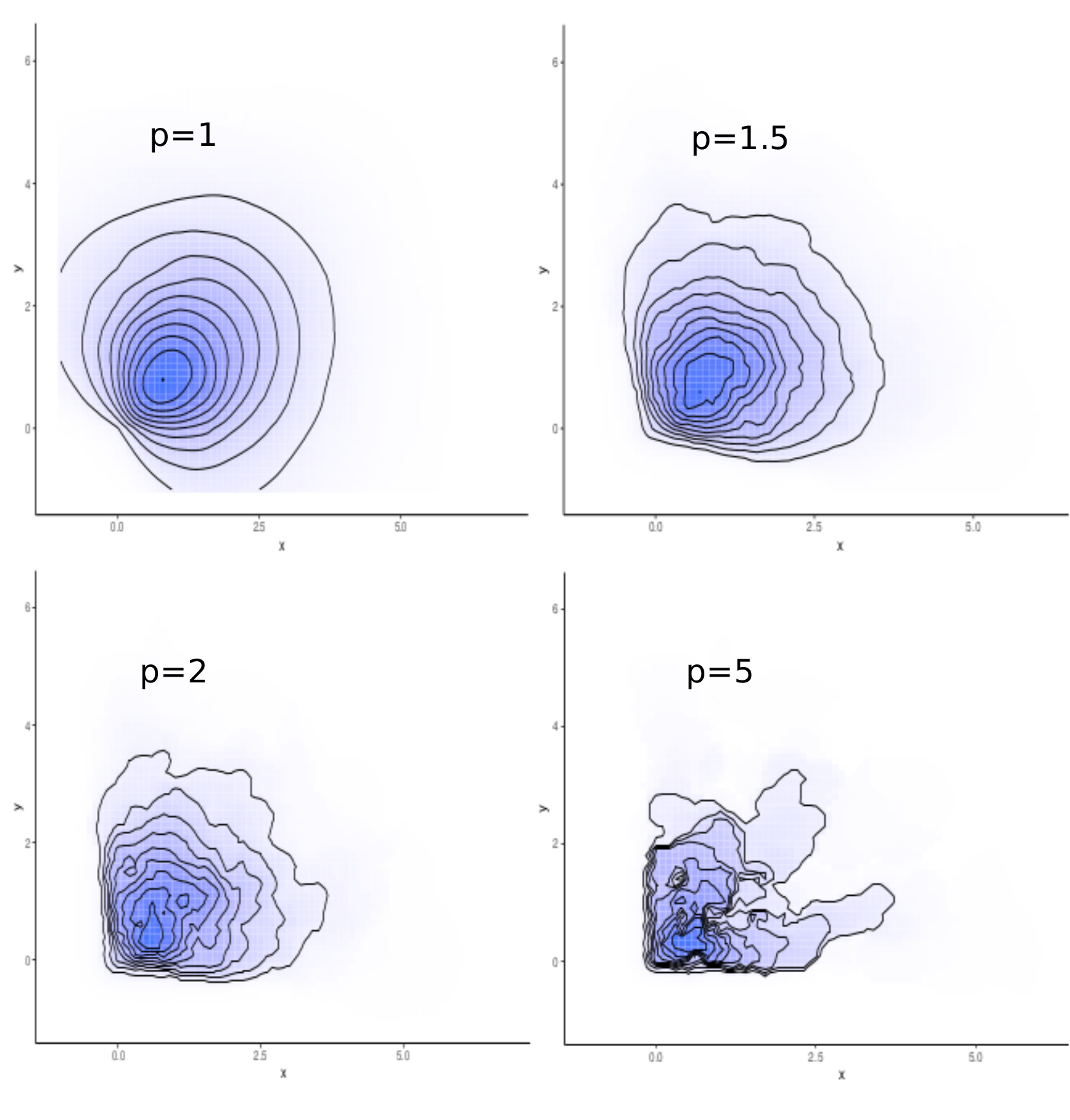}\label{exponencial} }
  \hfill
  \subfloat[]{\includegraphics[width=0.42\textwidth]{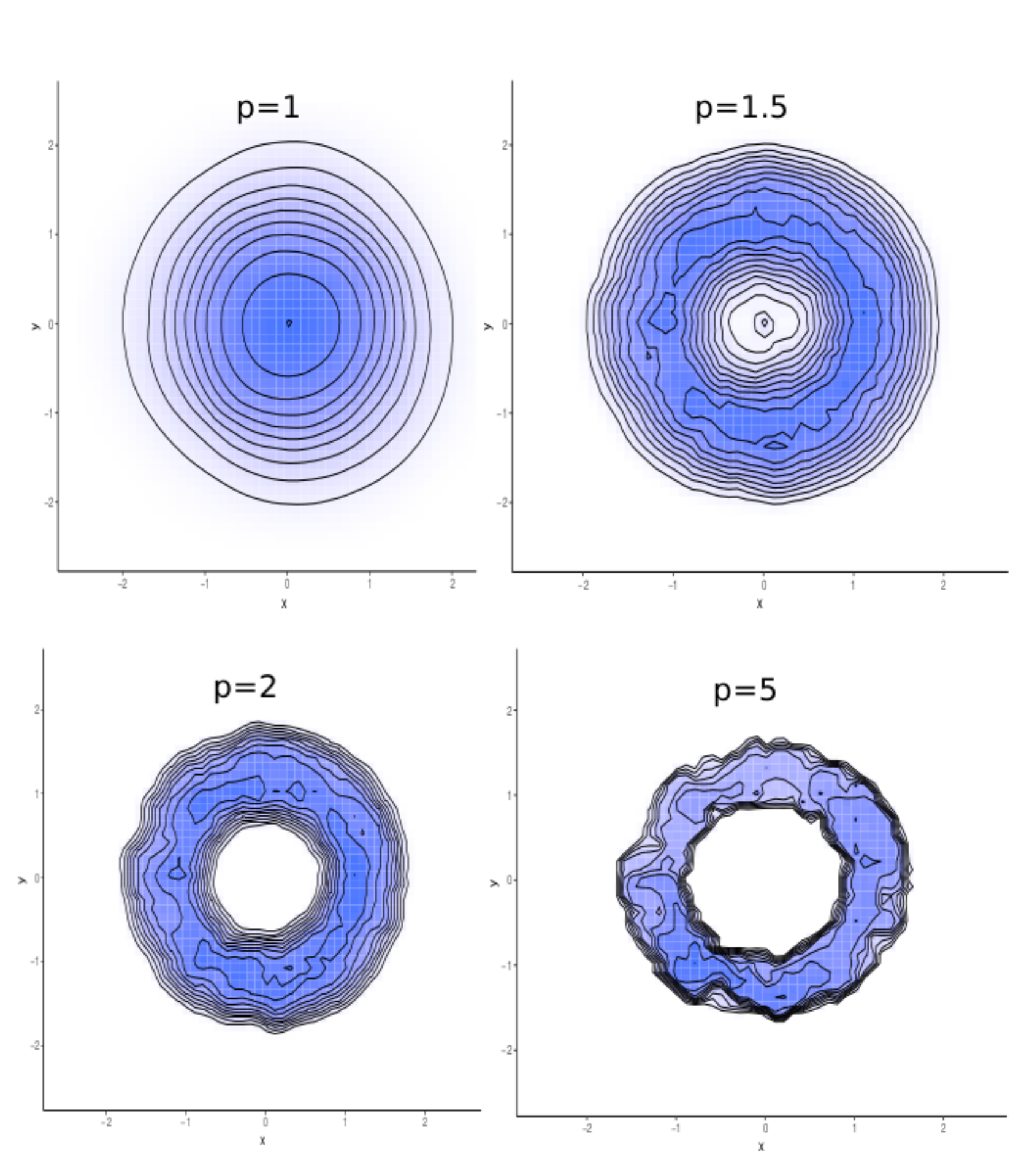}\label{uniforme} }
  \caption{(a) Level sets of a sample of $400$ exponentially distributed data-points.
(b) Level sets of a sample of $400$ uniformly distributed data-points on a ring.}
\end{figure}

Figures  \ref{exponencial} and \ref{uniforme} show how with increasing $p$ the level sets of $WLD$ adapt to the support of the distribution, even for non-convex level sets. This is not the case with  classical depths such as Tukey's or Liu's one.

\section{Classification with depth-depth-G and  $\textrm{WLD}_p$ }
\label{gg}
The idea to use a depth measure as dimension reduction tool to perform classification was considered in \cite{liu1990}. The method is called the depth-depth method, and is based on the depth-depth plots.
This method was improved in  \cite{li2012} and  \cite{cuesta2017}.  In these proposals,  after a reduction to dimension two (implementing the depth-depth method), a classifier is applied.
We will refer to this method as  depth-depth-G.

To explain the procedure, let us assume that we have two classes and a depth $D$. Then the depth-depth method assigns to each point $x$ a two dimensional vector, whose coordinates are its depth in each group. That is, $x \rightarrow \left(D_1(x), D_2(x) \right)$, where, for $i=1,2$, $D_i(x)$ is the depth of $x$ in the class $i$.
Finally, a classification method (for instance SVM, Random Forest, $k$-NN, neural networks) is applied, in dimension two.
This procedure can be easily generalized  to several classes,  or to product spaces, considering a depth on each component.
Moreover, more than one depth can be used for each space.
 
This method is a dimensionality reduction procedure and also a new and useful visualization technique, see \cite{ali2016}.

\subsection{Simulations}
\subsubsection*{Example 1: Interlocking rings (in $\mathbb{R}^2$)}

\begin{figure}[!ht] 
\centering
\subfloat{\includegraphics[width=76mm]{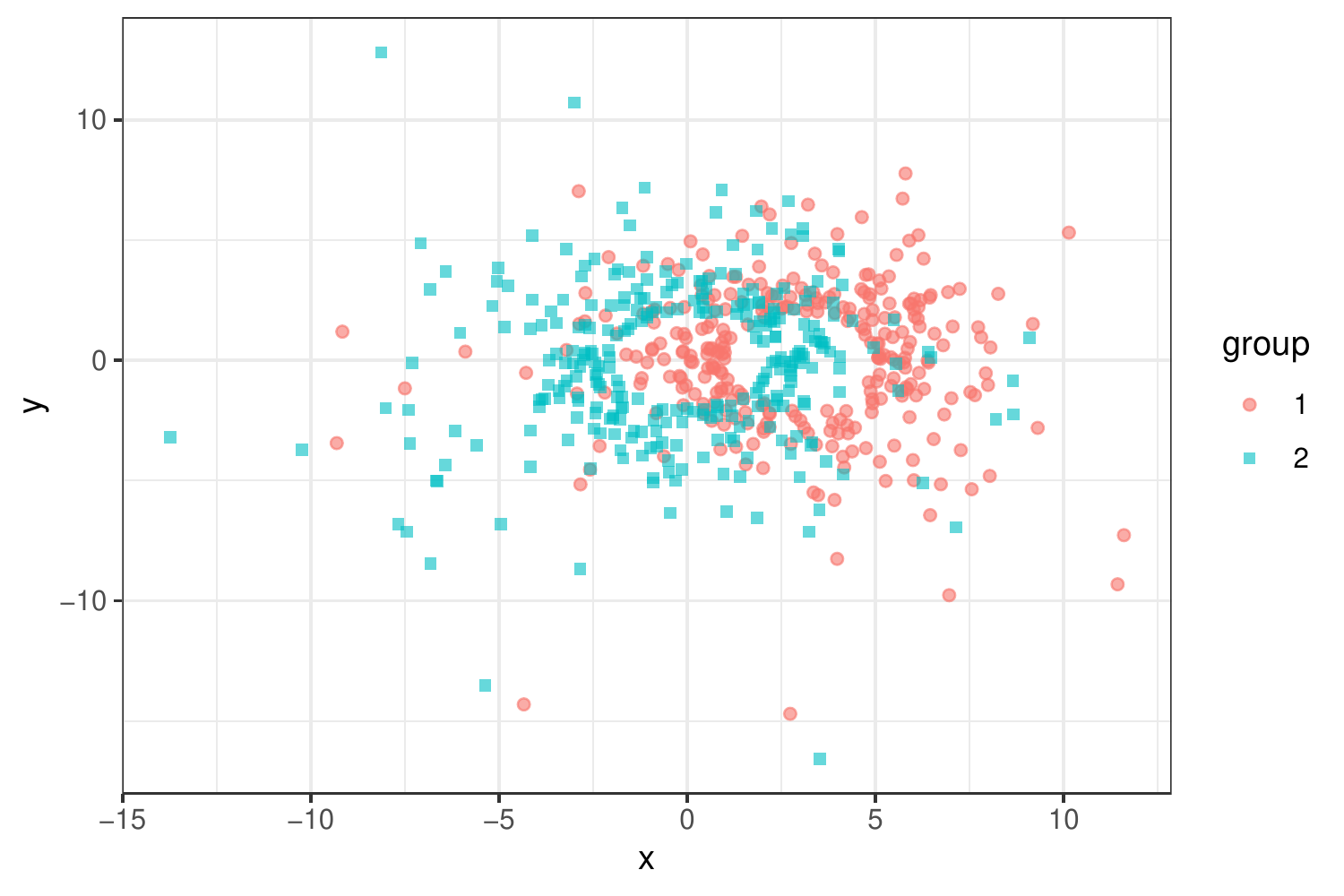}}
\caption{Sample points from the `interlocking rings' model } \label{anillosx}  
\end{figure}

In this example, we show how the value of $p$ captures different patterns that contribute to the classification.
We  consider as an example  the  `interlocking rings' model.
 Let $\mathbf{X}$ be a two-dimensional random vector, defined by $\textbf{X}= \left( (2+R) \cos \theta, (2+R) \sin \theta  \right)$, where $R$ and $\theta$ are independent random variables, with $R \sim \textrm{Exp}(1/2)$ and $\theta \sim \textrm{Unif} (0, 2 \pi)$.
We generated $300$ independent copies with distribution $\textbf{X}_1= \textbf{X}_1^{'}+ (3,0)$  (group 1) and $300$  from $\textbf{X}_2$ (group 2), where $\textbf{X}_1^{'}$ and $\textbf{X}_2$ are independent copies of $\textbf{X}$, (see Figure \ref{anillosx}).
We considered $70\%$ of the data for the training sample, and $30\%$ for validation.

\begin{figure}
    \begin{tabular}{p{0.28\textwidth}p{0.58\textwidth}}
\includegraphics[height=4cm,width=\hsize]{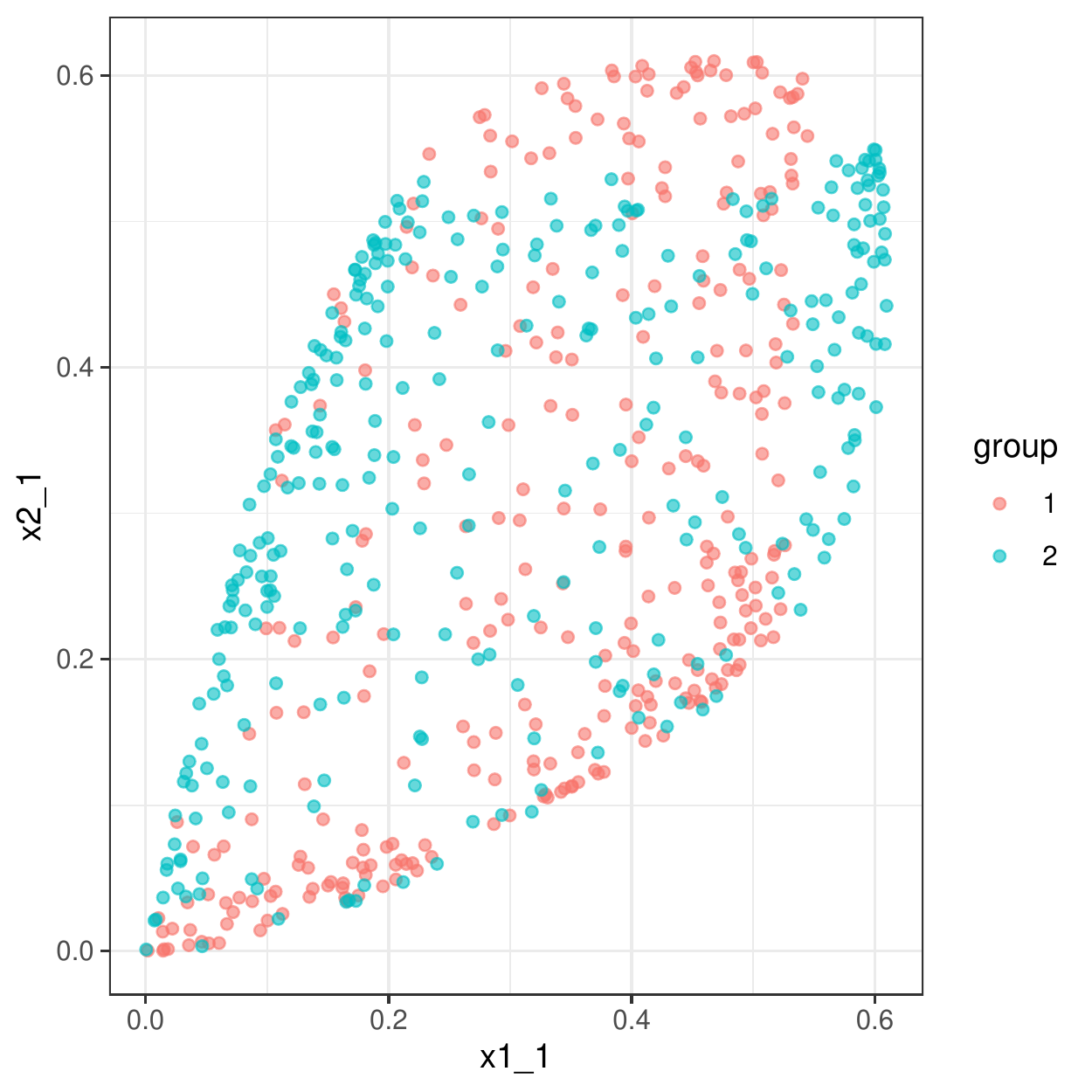}
&   \multirow{2}{\hsize}[37mm]{
    \includegraphics[height=6cm,width=\hsize]{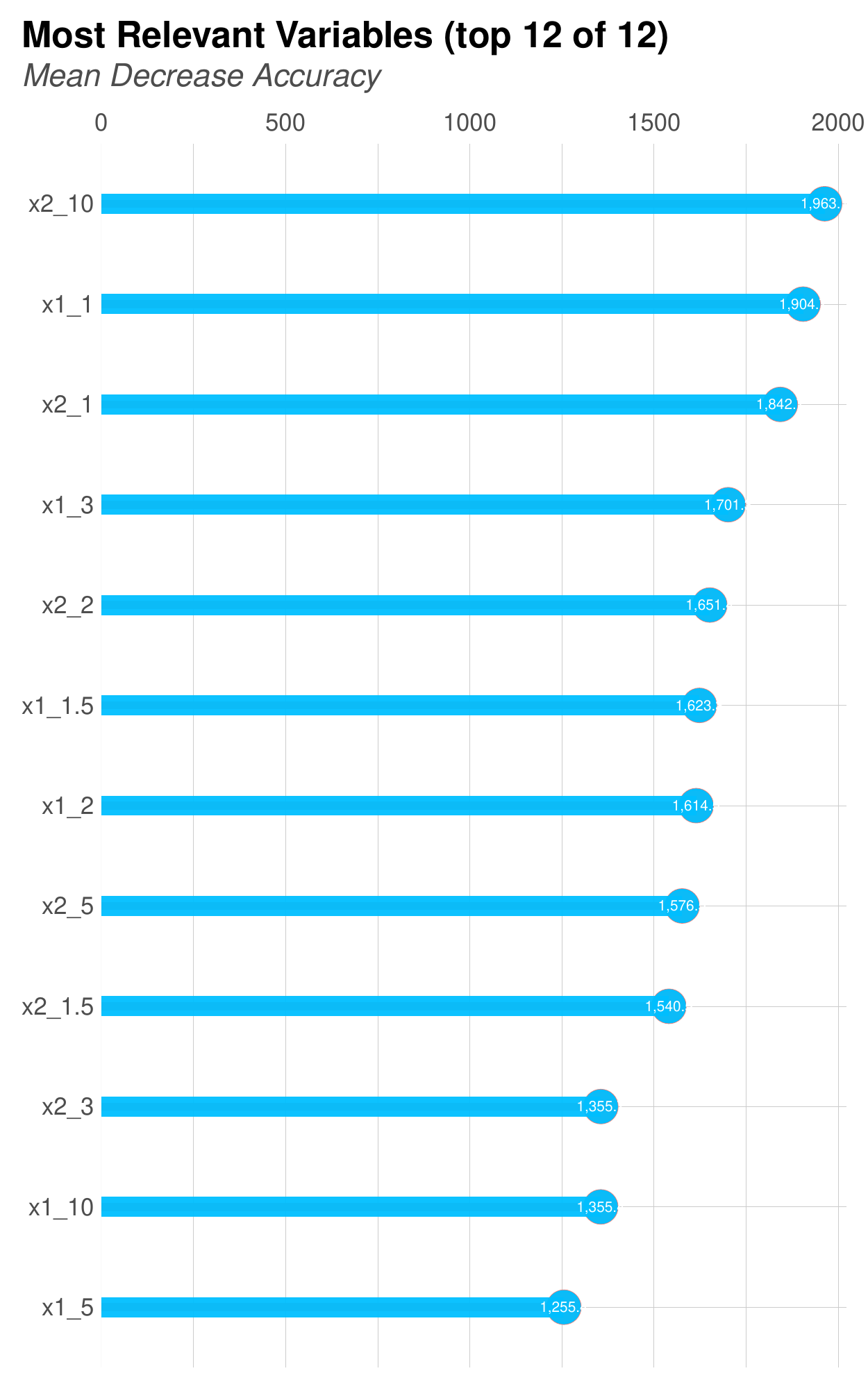}
}\\
\includegraphics[height=3cm,width=\hsize]{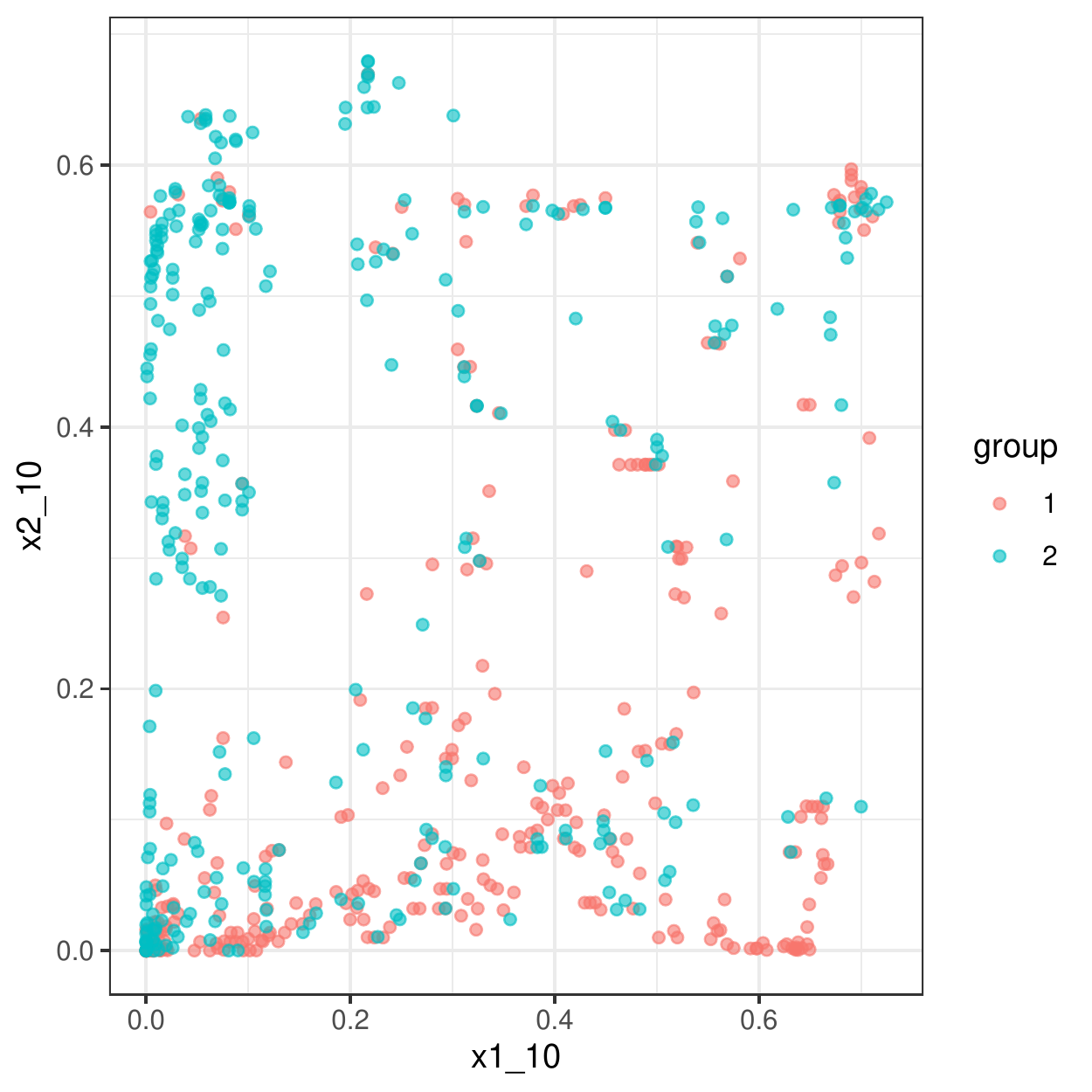}  &
    \end{tabular}
\caption{Top left Panel: Depth-depth plot for  $\textrm{WLD}_1$.
Botton left panel: Depth-depth plot for   $\textrm{WLD}_{10}$.
Right panel: Variable importance for RF.
}
\label{sss}
\end{figure}

We compared the performance of  1) Random Forest (RF), with the original sample; 2) depth-depth method with lens depth, and as the second step,  RF  in dimension $2$; 3) depth-depth-G with $\textrm{WLD}_p$ and $p= 1$ and $10$, and as the second step, RF in dimension four; 4) depth-depth-G with  $\textrm{WLD}_p$  and $p= 1,1.5,2,3,5,10$, and as the second step, RF in dimension $12$. Figure \ref{sss} shows the depth-depth plots  for $p=1$ and $p=10$.

The misclassification errors are shown in Table \ref{tanillos}.
 The main conclusion is that, depth-depth-G with  $\textrm{WLD}_p$ and $p= 1,1.5,2,3,5,10$ outperforms the other competitors.

\begin{table}[ht]
\caption{Misclassification errors}
\begin{center}	\footnotesize
\begin{tabular}{c|ccc}
   \toprule[0.4 mm]
    1) RF &\multicolumn{3}{c}{depth-depth with  $\textrm{WLD}_p$  and RF }\\
    \cline{2-4}
    &$2)\,\, p=1$ & $3)\,\, p=1,10$ &  $4)\,\,p= 1,1.5,2,3,5,10$ \\
     \cline{2-4}
    $0.28$ &$0.26$ & $0.22$ & $0.19$\\
  \toprule[0.4 mm]
\end{tabular}
\end{center}
\label{tanillos}
\end{table}

Let us consider case 4), whose classification error is the smallest.
We  studied the importance of the variables in the classification, using \textit{Mean Decrease Accuracy}, (see \cite{breiman2001}).
 Indeed,  they depend on the parameter $p$, see Figure \ref{sss}.
This also suggests that taking differing values of $p$ in $\textrm{WLD}_p$ leads to differing contributions to the classification.

\subsubsection*{Example 2: The cone of positive definite matrices}
 
Consider $\mathcal{C}_k$, the set of all $k\times k$ positive definite  matrices.
We will compare depth-depth-G using $\textrm{WLD}_p$  with RF,  on a binary classification problem on $\mathcal{C}_k$.
Given two matrices $A,B \in \mathcal{C}_k$, the geodesic distance is given by $d(A,B)= \Vert \ln  \left(A ^{-1/2} B A ^{-1/2}  \right) \Vert,$ where $\Vert \cdot \Vert$ is the Hilbert--Schmidt norm (i.e., $\|A\|=trace(AA^t)$ being $A^t$ the transpose of $A$), see \cite{moakher2005}. Also, there exists a unique geodesic joining $A$ and $B$ (see \cite{moakher2005}), given by $\gamma(s) \defeq  A ^{1/2} \left(A ^{-1/2} B A ^{-1/2}  \right) ^ s A ^{1/2}.$ 
To generate matrices with the Wishart distribution, we first choose a covariance matrix $\Sigma$ and a positive integer $m$. Second, we generate a sample of $m$ iid random variables with common distribution $N(0,\sigma)$.
Then the matrix $S \defeq \sum_{i=1}^m X_i X_i^\top$ has a Wishart distribution, i.e. $ S \sim \mathcal{W} (\Sigma, m)$ on  $\mathcal{C}_5$.

\begin{figure}[!ht]
\centering
\subfloat{\includegraphics[width=85mm]{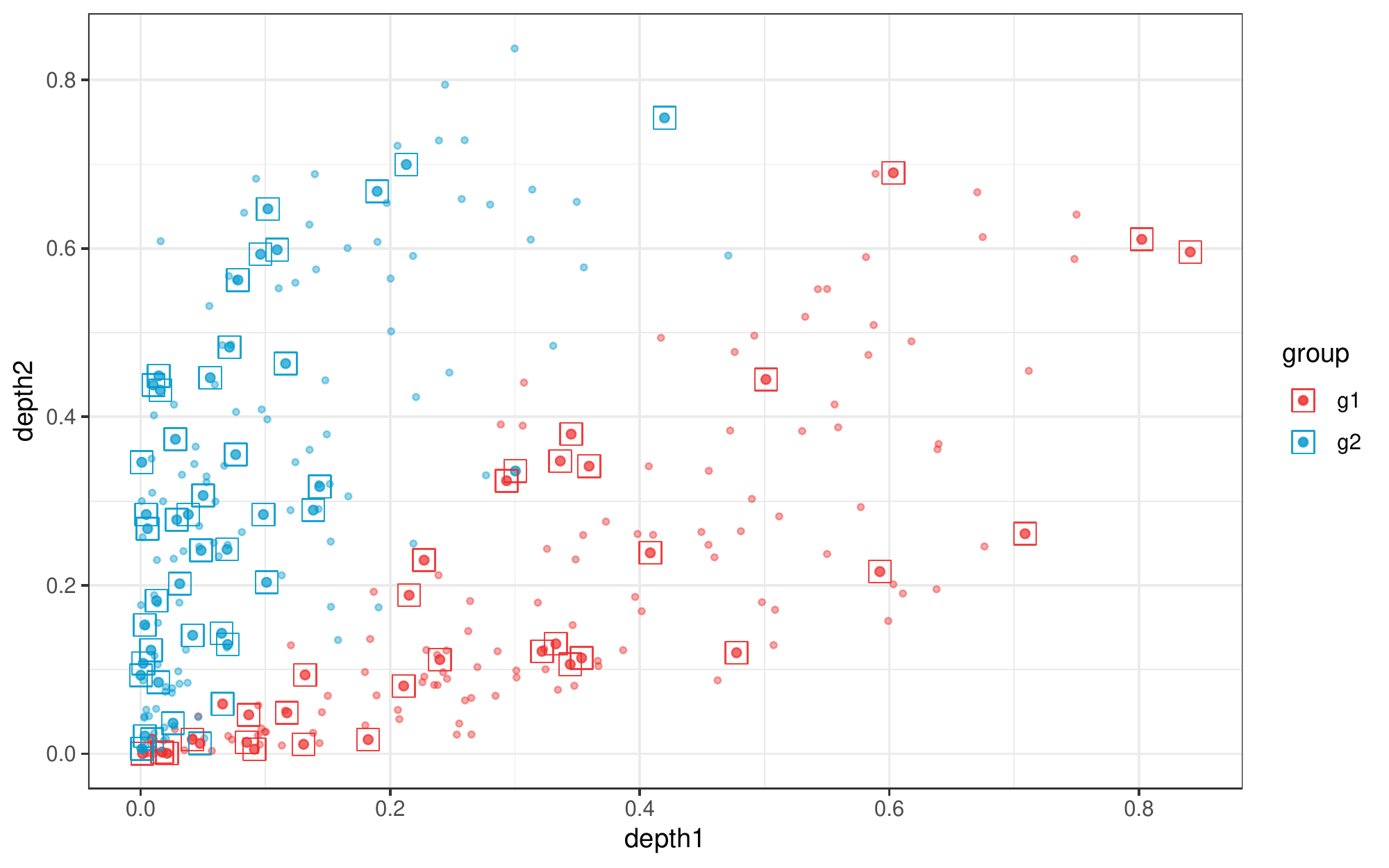}}
\caption{The result of the classification based on depth-depth-G with $\textrm{WLD}_2$ and RF, for $n=300$.
The colour of the circle shows the group to which each sample belongs.
The square shows the validation sample.
The colour of the square correspond to the group assigned.} \label{cara1}  
\end{figure}

In $\mathcal{C}_5$, we will consider two iid samples of size $n$, following  different Wishart distributions, each of them corresponding to a different group.
For the first group we choose $m_1=10$ and $\Sigma_1=I_5$.
For the second group we choose $m_2=10$ and   $\Sigma_2= 2I_5$.
The sample sizes  in our study are  $n=100,200,300$.
From each sample we choose $25\%$ for validation.
For the depth-depth method, we used first the lens depth, and as the second step RF.
Figure \ref{cara1} we display the depth-depth plot for $n=300$, with  $\textrm{WLD}_{2}$, and RF.
The misclassification error in that figure is $7\%$.
We performed a simulation study to compare depth-depth-G and  $\textrm{WLD}_p$  with RF for $p=1$,$1.5$,$2$,$5$, $k$-nearest neighbour ($k$-NN) classifier for $k=[n/10]$.
The results are displayed in Table \ref{dd}.
The whole procedure was replicated 100 times.
 The conclusion is that  depth-depth-G and  $\textrm{WLD}_p$ with RF outperforms $k$-NN in all cases, and the best value for the parameter $p$ is $p=2$.
\begin{table}
	\begin{center}
	\footnotesize
		\begin{tabular}{c|cccc|c}
 \toprule[0.4 mm]
\multirow{3}{*}{Sample size} & \multicolumn{4}{c|}{depth-depth-G and  $\textrm{WLD}_p$ , with RF} & %
   \multirow{3}{*}{$k$-n.n}\\
\cline{2-5} 
 & $p=1$ & $p=1.5$ & $p=2$& $p=5$  & \\
 \toprule[0.1 mm]
$n$=100& 0.20 & 0.14  & 0.11 & 0.18 & 0.24 \\
 \toprule[0.1 mm]
$n$=200& 0.12 & 0.08  & 0.09 &  0.16 & 0.20 \\
 \toprule[0.1 mm]
$n$=300 & 0.10 &  0.08 & 0.09&  0.18  & 0.21\\
 \toprule[0.4 mm]
\end{tabular}
\caption{Misclassification error for the positive definite matrices} \label{dd}
\end{center}
\end{table}

\subsection{An analysis of some real-life data}

\subsubsection{The manifold of phylogenetic trees}

\begin{figure}[!ht] 
\centering
\subfloat{\includegraphics[width=100mm]{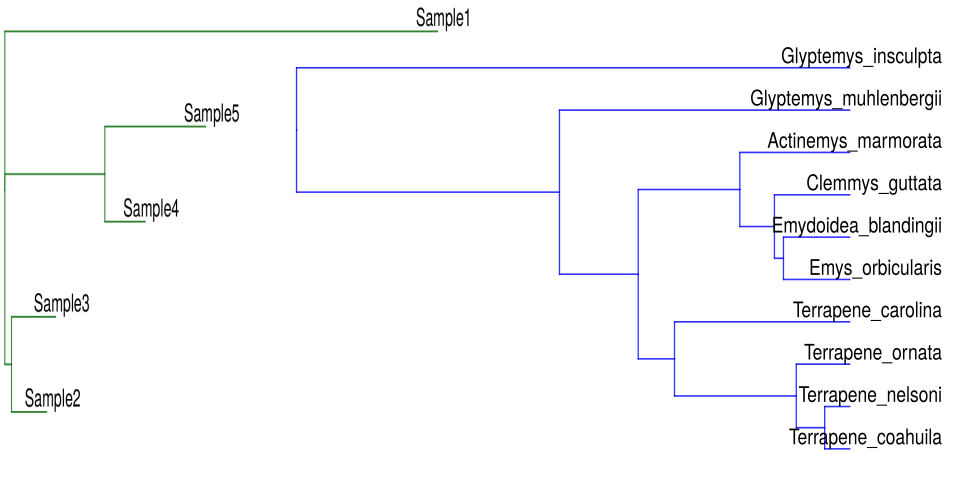}}
\caption{Left tree: A sample from the year 2001.
Right  tree: A sample from the gen ODC.
} \label{arbb}  
\end{figure}

In this section,  we focus on two examples arising from genetics.
In both examples, the information contained in a genetic chain is approximated, represented, and analysed by means of a sample of low-complexity phylogenetic trees.
 These trees are used in general to represent ancestral histories.
For details on how to build these trees and how to endow them with a Riemannian metric we refer to  \cite{billera2001}.
This manifold is complete and separable.
An example of two of these trees is shown in Figure \ref{arbb}.

\subsubsection*{Example 3: Temporal analysis of influenza}

As the first example, we use depth-depth-G on the manifold of phylogenetic trees, to analyse the evolution over time of the influenza virus.
 Modelling the genomic evolution of these viruses is an important issue ( see \cite{smith2004b}).
They mutate and evolve constantly, becoming resilient to different drugs.
Some examples are the Avian flu, SARS, and H2N2, among others.
To be able to predict new genetic patterns is crucial to developing vaccines.
 Several statistical methodologies have been presented in the literature to understand this phenomenon, for example, in  \cite{solovyov2009}, clusters of chains of genes are built to understand the pandemic of 2009.

We will focus on influenza of type A, H3N2, in particular, the subtype hemagglutinin (HA), since this subtype presents the largest variability in its genomic evolution and a strong resistance to standard vaccines, see \cite{altman2006}.
 The data used here   can be found in the GI-SAID $\textrm{EpiFlu}^{TM}$ database \footnote{\url{www.gisaid.org}},
which consists of the genomic data of $1089$ sequences from 1993 to 2017 collected in New York.
 A pre-processing of the data is required, details of which can be found in \cite{zairis2016, monod2018}.
 
The data-base was obtained from the GitHub repository \url{https://github.com/antheamonod/FluPCA}.
Figure \ref{arbb} (left tree) shows a tree from the sample of 2001.

\begin{figure}[!ht] 
\centering
\subfloat{\includegraphics[width=140mm]{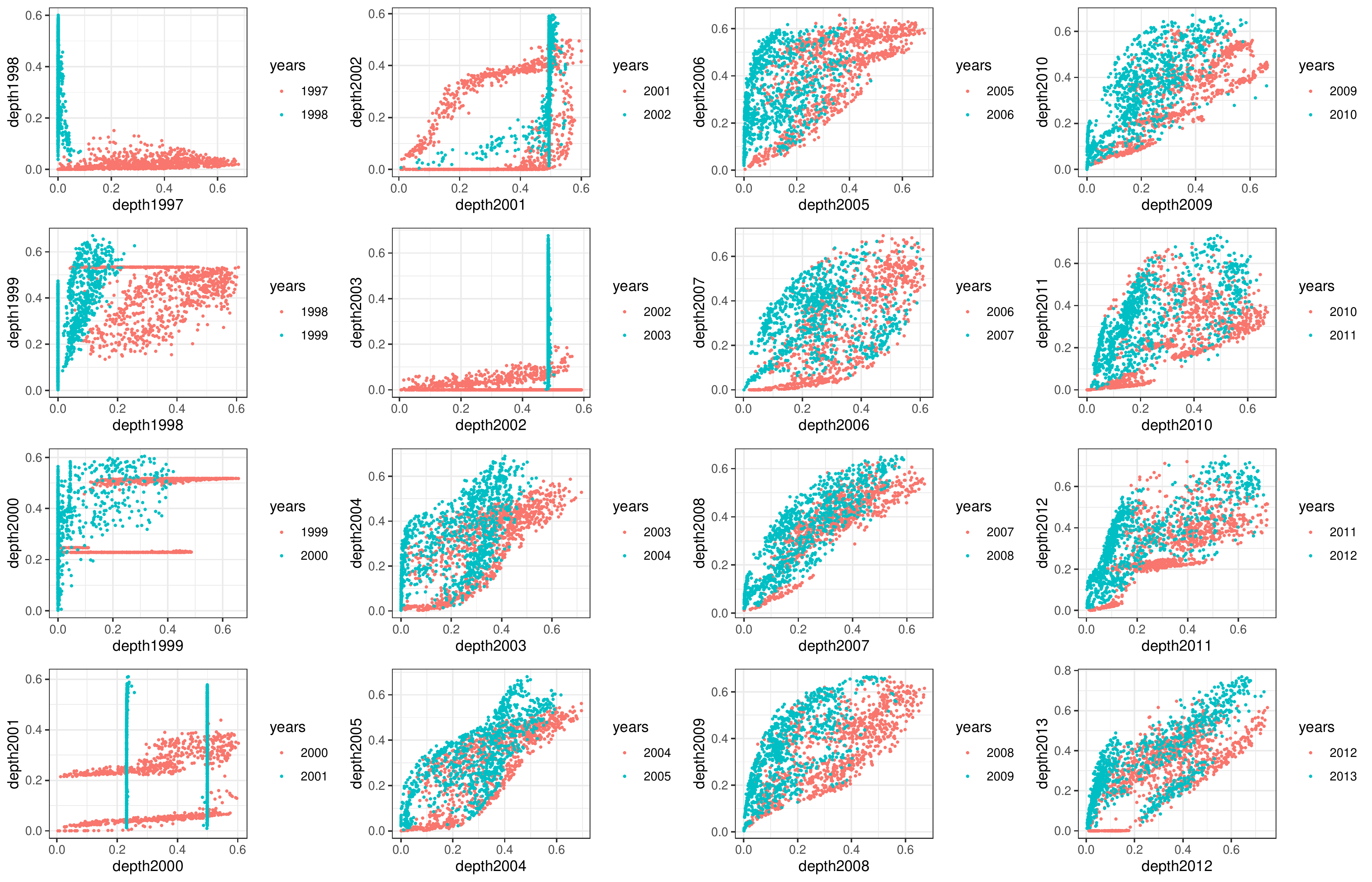}}
\caption{Depth-depth plots with lens depth for the sample of phylogenetic trees.} \label{genetica1}  
\end{figure}

We will use as a dimensionality reduction tool the depth-depth method between consecutive years, using LD with the geodesic distance.
This is different from the PCA used in \cite{monod2018,nye2017}.
It allows a better visualization of the evolution of the different genetic orders, and the new kinds of  strains, see Figure \ref{genetica1}.

 We may observe  a clear difference between the strains of the virus from one year to the other, until 2003.
Moreover, it seems that in some years (like 1999 or 2001), two different strains coexist.
From 2003, the strains became more similar from one year to the other.
It worth mentioning that the pandemic in 2009 is not detected by the diagrams because it involved H1N1 and we are analysing H3N2.

\subsubsection*{Example 4: Evolutionary trees for species of turtles}

\begin{figure}[!ht] 
\centering
\subfloat{\includegraphics[width=66mm]{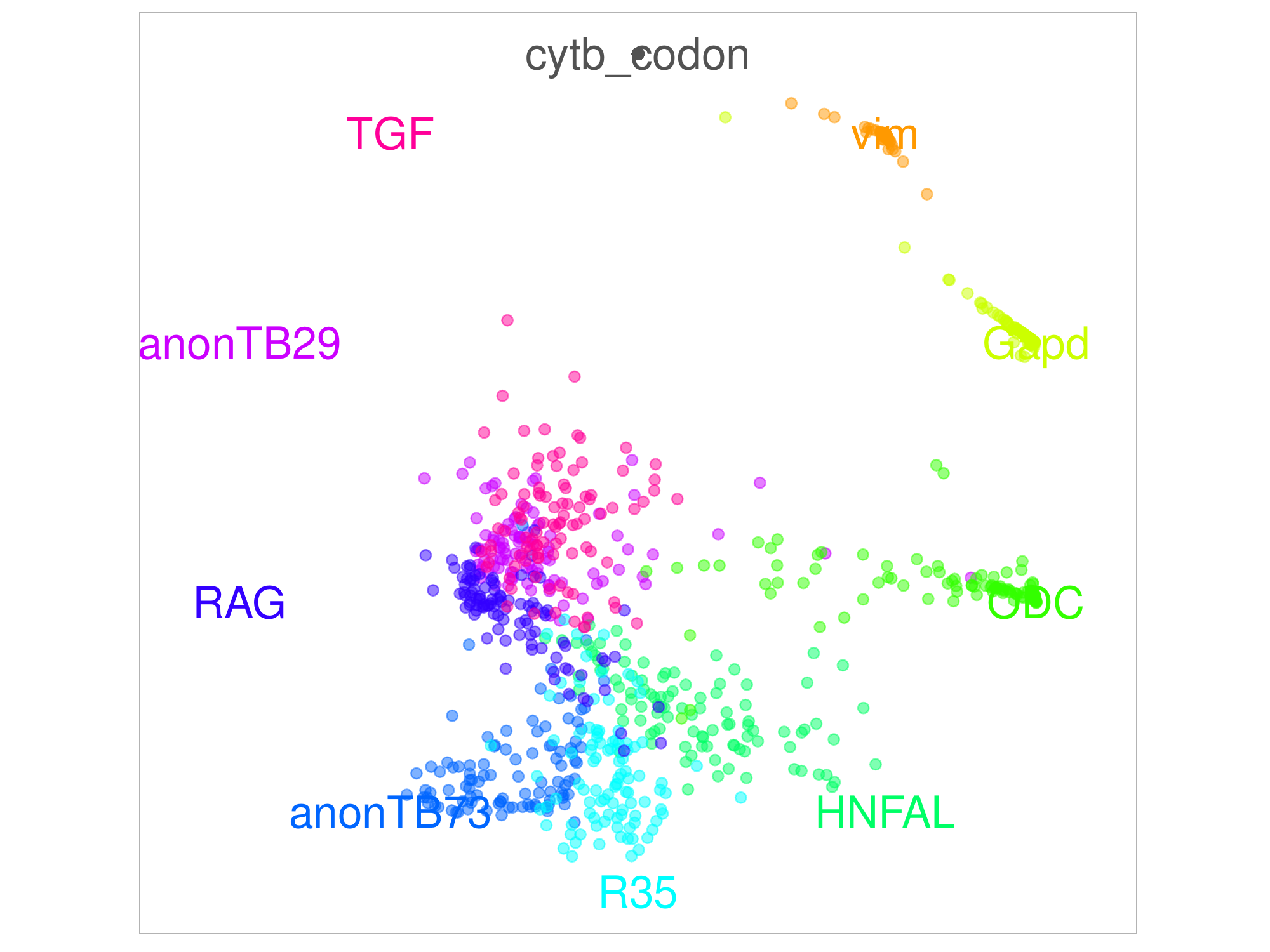}}
\subfloat{\includegraphics[width=72mm]{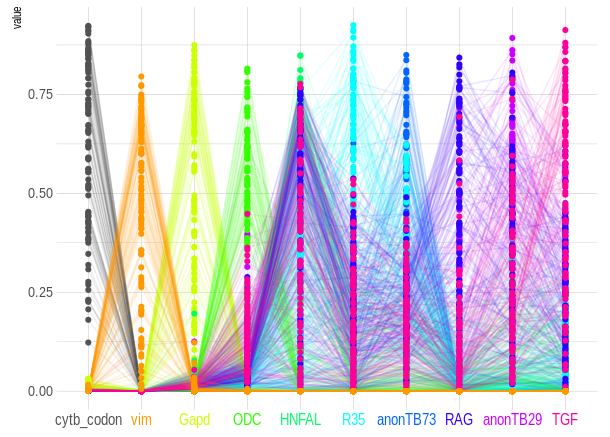}}
\caption{Depth for the sample of phylogenetic trees for species of turtles.
Left panel: Weighted circular depth diagram for relative depths.
Right panel: Parallel coordinate plot.} \label{tortugas}  
\end{figure}

We will apply a similar approach  to the evolutionary trees of the genes of species of turtles.
The data-set consists of 1000 evolutionary trees. Each tree corresponds to one of ten possible genes (nine of them are of nuclear type, and the other is of mitochondrial type).
There are 100 trees for each gene.
The aim is to visualize the difference between the genes using lens depth.
This problem is particularly difficult due to the large number of nodes on the evolutionary trees, (see \cite{wiens2010}).
 A pre-processing of the data is performed (see  \cite{willis2018} for details).
Figure \ref{arbb} (right) shows a tree of the species ODC.
Using the depth-depth method with the lens depth on each gene, we have been able to find a clear difference between  the mitochondrial type gene and the other nine nuclear genes.
We have also been able  to catch the differences between some nuclear genes, which the PCA used in \cite{wiens2010} with two components was not able to capture.  Figure \ref{tortugas} shows two different representations of the depths obtained (using LD): at the left, a circular diagram of the relative depths.
Roughly speaking, the position of each point on the circle is based on the compositional-data depth of a tree on a given group, w.r.t. to its own group relative to the other nine groups  (see \cite{hoffman1999}  for details of how these diagrams are built).
 On the right we show the parallel coordinate plot (we used the function ggparcoord of the R package ggally), that represents the absolute depths of the trees.  Notice that  the mitochondrial gene cytb-codon  can be clearly differentiated from the other, nuclear, genes.
Some genes present very different depths when compared to the others, for instance  the genes vim, Gapd and ODC.
However, some of them (such as TGF) present very high depths on their own gene, but also on the other genes, as can be seen in the parallel coordinate plot.
 This made the visualization problem much more involved.

\subsubsection{A signal recognition problem}

A widely studied problem in signal recognition is how to automatically identify the nationality of a person when the available data is a sample of recorded words pronounced by a group of people, (see for instance \cite{dhanalakshmi2009}).
We used a data-base consisting of 329 people, from 6 different nationalities, where each of them pronounces a word in English.

\begin{figure}[!ht] 
\centering
\subfloat{\includegraphics[width=110mm]{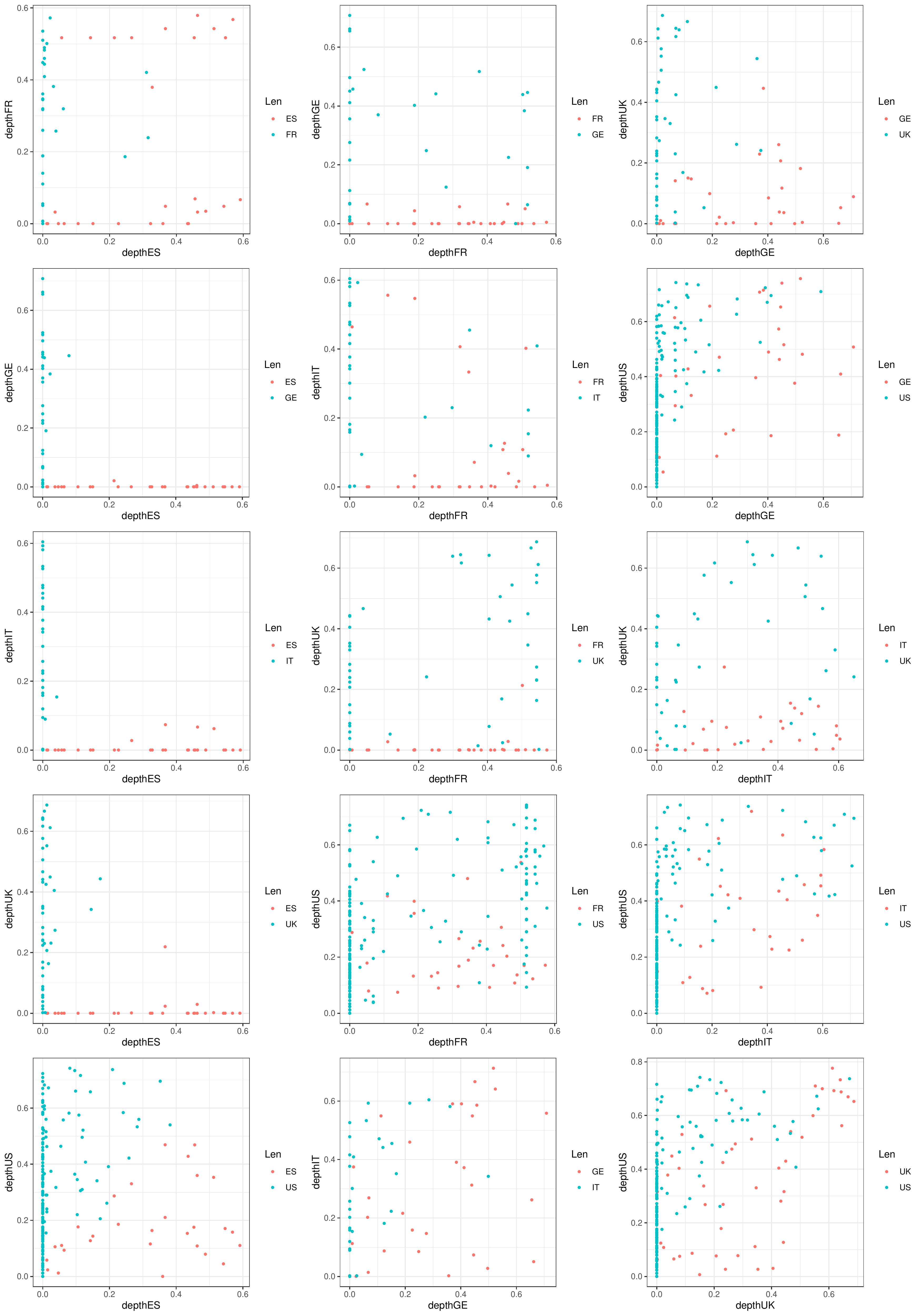}}
\caption{Depth-depth plots for  $\textrm{WLD}_5$, for all pairs of nationalities considered.} \label{genetica2}  
\end{figure}

The data-base can be downloaded from \url{https://archive.ics.uci.edu/} (see \cite{fok2020}).
Each signal was preprocessed using the 
MFCC (Mel-Frequency Cepstral Coefficients) method, see \cite{pedersen2008} and \cite{ma2015}.
The countries considered were Spain, France, Germany, the USA, Italy and England, with sample sizes $29$, $30$, $30$, $165$, $30$ and $45$, respectively. There is a reasonable balance with respect to the gender of the participants. As can be seen in Figure \ref{genetica2}, a good performance is obtained by the depth-depth method with  $\textrm{WLD}_5$.
For instance, there seems to be a clear separation between the Spanish and other nationalities, such as German or North American.
In other cases, the pronunciation is quite similar, for instance between Germans and North Americans.
One can guess that this is due to the common evolutionary root of both languages.

\section{Concluding remarks}

 Studying  depths on general metric spaces allows tackling problems where the structure of the data is neither finite dimensional nor functional. Some important examples, among others, are the phylogenetic trees given by the  evolution of the history of genes, or data belonging to an unknown manifold. For separable and complete metric spaces, we have studied the main properties of the lens depth (LD), proved the consistency of the empirical version with the population one, and have provided almost parametric convergence rates. For Riemannian manifolds, we have provided an extension of LD, called $\textrm{WLD}_p$, which is  more flexible than LD, and is able to also catch information about the geometry of the underlying distribution.
We have also proved consistency results for the plug-in estimator of $\textrm{WLD}_p$. For most of  the aforementioned problems, we showed, by using $\textrm{WLD}_p$, that it is crucial that the depth take into account not only the geometric structure of the data but also the underlying distribution. We use $LD$ and $WLD_p$ in supervised classification problem, putting in action the depth-depth method introduced in \cite{liu1990}.  These classifications are performed on simulated and real life data examples. We obtain  better performance than some competitors and interesting results on the real life data examples.
 
 
\section*{Appendix}

\textit{Proof of Proposition \ref{prop1}}

	\begin{itemize} 		\item[(a)]
		
		Let $x_1,x_2\in M$, and write $R=d(x_1, x_2)$.
Choose $K$ large enough to ensure  $x_1,x_2 \in B(y, K/3)$.
If $d(x,y)>K$, it follows that $d(x_1,x)>2K/3$.
Then $x \notin B(x_1, R) \cap B(x_2, R)$.
If $X_1$ and $X_2$ are independent copies of $X$, then 
$\textrm{LD}(x)=\mathbb{P}  \left( x \in B(X_1, R) \cap B(X_2, R) \right) \leq 2\mathbb{P} (X \notin B(y, K/3))$ 	and $\mathbb{P} (X \notin B(y, K/3))\to 0 $ as $K\to \infty$.
\item[(b)]
		Let us assume that  (b) does not hold. Then there exists $\Omega_1\subset \Omega$ with $\mathbb{P} (\Omega_1)=\delta>0$ such that for all $\omega\in \Omega_1$ there exists $\gamma=\gamma(\omega)>0$, and a sequence $K_n\to \infty$ such that $\sup_{\{x:d(x,y)>K_n\}} \widehat{\textrm{LD}}_n(x)>\gamma>0.$
		Then for all $\omega\in \Omega_1$ and for all $K$, there exists $x=x(K)$ such that $d(x,y)>K$  and $\widehat{\textrm{LD}}_n(x)>\gamma$.
Let $K_0$ be large enough to ensure $\mathbb{P}(d(y,X)>K_0/4)< \delta/(2n)$.
Since $\Omega_1 \subset \left\{ \omega : \text{ there exists } i \text{ such that } d(y,X_i(\omega))>K_0/4\right\},$		it follows that  $\delta \leq n \mathbb{P} \left(d(y,X)>K_0/4 \right)< \delta/2$, which is a contradiction.\QEDB
	\end{itemize} 

\textit{Proof of Proposition \ref{prop2}}
	\begin{itemize}
		\item[(a)]
$\vert  {\textrm{LD}}(x) -  \textrm{LD}(y)  \vert  =  \int_{M^2} |\mathbb{I}_{A(z,t)}(x)-\mathbb{I}_{A(z,t)}(y)|P_X(dz)P_X(dt).$
		By definition, $1=|\mathbb{I}_{A(z,t)}(x)-\mathbb{I}_{A(z,t)}(y)|$ if, and only if, either both $x\in A(z,t)$ and $y\notin A(z,t)$, or $x\notin A(z,t)$ and $y\in A(z,t)$.
Set  $d(x,y)=\epsilon$.
In the first case, $x\in [\partial B(z,d(z,t))]^\epsilon\cup [ \partial B(t,d(z,t))]^\epsilon=:\mathcal{E}_{z,t}^\epsilon$ and in the second case $y\in \mathcal{E}_{z,t}^\epsilon$
Then
		\begin{equation*}
			\int_{M^2} |\mathbb{I}_{A(z,t)}(x)-\mathbb{I}_{A(z,t)}(y)|P_X(dz)P_X(dt)\leq 
			\int_{M^2} \mathbb{I}_{\mathcal{E}_{z,t}^\epsilon}(x)P_X(dz) P_X(dt).
		\end{equation*}
		For any $z,t$, and $x$, $\mathbb{I}_{\mathcal{E}_{z,t}^\epsilon}(x)\to \mathbb{I}_{\mathcal{E}_{z,t}^0}(x)$ as $\epsilon\to 0$, so 
		 by the dominated convergence theorem, the  integral converges to $2\mathbb{P}(x\in \partial B(X_1,d(X_1,X_2))).$ Lastly,  \eqref{lcont}, implies that $\mathbb{P}(x\in \partial B(X_1,d(X_1,X_2)))=0$.

		\item[(b)]
		The proof is a straightforward consequence of the fact that the kernel of the U-statistic of order two  $\widehat{\textrm{LD}}_n(x)$ is bounded between $0$ and $1$, together with Hoeffding's inequality for U-statistics, (see \cite{serfling1980}, p.~201). This implies that for all $\delta>0$ and $n>2$,
		\begin{equation}
			\label{ofi}
			\mathbb{P} \left(  \beta_n\vert \widehat{\textrm{LD}}_n(x) - {\textrm{LD}}(x) \vert > \delta \right) \leq 2 \exp\Big(-\frac{n}{\beta_n^2}C\delta^2\Big), 
		\end{equation}
		with $C>0$. Therefore, we may conclude applying the Borel--Cantelli Lemma. \QEDB
	\end{itemize}
\textit{Proof of Lemma \ref{lema1}}

Let us define, for $\epsilon>0$, $\Omega_1=  \{ \omega : \sup_{x_1,x_2 \in M, \rho_1(x_1,x_2) \geq b}
   \vert \frac{L(x_1,x_2; \mathcal{X}_n)}{ n^{(1-p)/d}\rho_p(x_1,x_2)} -C(d,p)   \vert > \epsilon \}, $
and $B= \overline{B_1(x_1,b)^c}$.
 Given $\delta>0$, let $0<b< \min \{ \delta/6, \rho_p(x_1,x_2) \} $.
By the triangle inequality,  
\begin{align*}
	\mathbb{P} \left( d_{H} \left[B_p(x_1, \rho_p(x_1,x_2)),  A^1_{p,n}(x_1, x_2) \right] > \delta  \right)& \\
	&\hspace{-6cm}\leq  \mathbb{P} (d_{H} \left[B_p(x_1, \rho_p(x_1,x_2)), B_p(x_1, \rho_p(x_1,x_2)) \cap B \right] > \delta/3)\\
	&\hspace{-6 cm}+ \mathbb{P}  \left( d_{H} \left[B_p(x_1, \rho_p(x_1,x_2)) \cap B, A^1_{p,n}(x_1, x_2) \cap B \right] > \delta/3  \right) + \\
	& \hspace{-6 cm} + \mathbb{P}  \left( d_{H} \left[A^1_{p,n}(x_1, x_2) \cap B, A^1_{p,n}(x_1, x_2) \right] > \delta/3 \right)
	= J_1 +J_2 +J_3.
\end{align*}
From $d_{H} \left[ B_p(x_1, \rho_p(x_1,x_2)) \cap B, B_p(x_1, \rho_p(x_1,x_2)) \right] \leq \textrm{diam} \left( B^c \right) <\delta/3$, it follows that  $J_1=0$.
In the same way, $J_3=0$.
Now $J_2$ is bounded from above by  
\begin{align*}
	\mathbb{P}  \left( d_{H} \left [B_p(x_1, \rho_p(x_1,x_2))\cap B, A^1_{p,n}(x_1, x_2) \cap B \right]  > \delta/3 \big \vert \Omega^c_1 \right) + \mathbb{P} (\Omega_1)= I_1 +I_2.
\end{align*}
By  Theorem \ref{teo1}, $I_2 \leq \exp \left(- \theta_0 n^{1/(d+2p)}\right)$.
Let us bound $I_1$.
By Theorem 1, for all $\omega\notin \Omega_1$,
$$n^{\alpha}\rho_p(x_1,x_2)(C(d,p)-\epsilon)\leq L(x_1,x_2;\mathcal{X}_n)\leq n^\alpha\rho_p(x_1,x_2)(C(d,p)+\epsilon),$$
$$ \mbox{ and } \ \ n^{\alpha}\rho_p(t,x_1)(C(d,p)-\epsilon)\leq L(x_1,t;\mathcal{X}_n)\leq n^\alpha\rho_p(t,x_1)(C(d,p)+\epsilon),$$
where $\alpha={(1-p)/d}$.
Let us define $C^+_-=(C(d,p)  + \epsilon)/(C(d,p)  -\epsilon)$ and $C^-_+= (C(d,p)  -\epsilon)/(C(d,p)  + \epsilon)$, then for all $t \in B$ 
\begin{equation}\label{eq00}
	C^-_+  \frac{\rho_p(t,x_1)}{ \rho_p(x_1,x_2)} \leq \frac{L(t,x_1; \mathcal{X}_n)}{L(x_1,x_2; \mathcal{X}_n)} \leq C^+_-\frac{\rho_p(t,x_1)}{ \rho_p(x_1,x_2)}.
\end{equation}

If we define $D^{+} \defeq B_p \left( x_1, C^+_-\rho_p(x_1,x_2) \right)$ and $D^{-} \defeq B_p \left( x_1, C_+^- \rho_p(x_1,x_2) \right)$ then from \eqref{eq00},
$D^{-} \cap B \subset A^1_{p,n}(x_1, x_2) \cap B \subset  D^{+}  \cap  B.$ Since $C^-_+<1$ and $C^+_->1$, $D^{-}  \cap B \subset B_p \big( x_1, \rho_p(x_1,x_2)\big)  \cap B  \subset D^{+}  \cap  B.$
Further, conditioned on $\omega \notin \Omega_1$, for $n$ large enough,
$$d_{H} \left [B_p(x_1, \rho_p(x_1,x_2))\cap B, A^1_{p,n}(x_1, x_2) \cap B \right]\leq d_{H} \left[ D^{-}\cap B,  D^{+}\cap B \right].$$ 
Observe that $ d_{H} \left[ D^{-}\cap B,  D^{+}\cap B \right]\to 0$ as $\epsilon\to 0$. So that, for $\epsilon$ small enough and $n$ large enough, $I_1=0$. \QEDB

\textit{Proof of Theorem \ref{conswld}} 

Let $\delta> 0$. $\mathbb{P}  (   \vert \widehat{\textrm{WLD}}_{p,n}(x)-\textrm{WLD}_p(x) \vert > 2\delta/\beta_n )$ is bounded from above by 
 	\begin{align*}
 		&  \mathbb{P}  \left(  \left  \vert  \binom{n}{2}^{-1} \sum_{1 \leq i_1 < i_2 \leq n} \I_{A_{p,n}(X_{i_1} X_{i_2})}(x) - \I_{A_{p}(X_{i_1} X_{i_2})}(x)  \right \vert > \delta/\beta_n  \right)  + \\
 		&  \mathbb{P}  \left(  \left  \vert  \binom{n}{2}^{-1} \sum_{1 \leq i_1 < i_2 \leq n} \I_{A_{p}(X_{i_1} X_{i_2})}(x)- E \left(\I_{A_p(X_1,X_2)}(x) \right)  \right \vert > \delta/\beta_n  \right) = I + II.
 	\end{align*}
 	
The term II   is bounded from above by $\exp(-Cn/\beta_n^2)$, arguing similarly as was done to prove  \eqref{ofi} (with $g_p$ instead of $g_1$), $C>0$ being a constant.  By Chebyshev inequality,	I is bounded from above by  $(\beta_n/\delta)\mathbb{E}  \left[ \I_{A_{p,n}(X_{1} X_{2}) \triangle A_{p}(X_{1} X_{2})}(x) \right] =(\beta_n/\delta)\mathbb{P} (x \in A_{p,n}(X_{1} X_{2}) \triangle A_{p}(X_{1} X_{2}) ).$
 	\begin{multline} 
 	\mathbb{P} \left(x \in A_{p,n}(X_{1}, X_{2}) \triangle A_{p}(X_{1}, X_{2})  \right) =\\
 		\int_{M^2} \mathbb{P}\big(x\in A_{p,n}(x_1,x_2)\triangle A_p(x_1,x_2)\big)P_X(dx_1) P_X(dx_2).
 	\end{multline}
 	Let us bound $\mathbb{P} (x\in A_{p,n}(x_1,x_2)\triangle A_p(x_1,x_2))\leq \mathbb{P} (x\in B(\partial A_p(x_1,x_2),\gamma),$
 	where $\gamma=d_H(A_{p,n}(x_1,x_2),A_{p}(x_1,x_2))$.
 Set $d_H=d_H(A_{p,n}(x_1,x_2),A_{p}(x_1,x_2)),$
 Let $P_{d_H}$ denotes its law on $\mathbb{R}^+$. Then,
 	\begin{multline*}
 		\int_{M^2} \mathbb{P}\big(x\in A_{p,n}(x_1,x_2)\triangle A_p(x_1,x_2)\big)P_X(dx_1)P_X(dx_2)=\\
 		\int_{M^2}\int_0^{+\infty}\mathbb{P}(x\in B_1(\partial A_p(x_1,x_2),\gamma))P_{d_H}(d\gamma)P_X(dx_1) P_X(dx_2)=\\
 		\int_0^{d_{\rho_1}(M)} \Big(\int_{M^2} P(x\in B_1(\partial A_p(x_1,x_2),\gamma))P_X(dx_1)P_X(dx_2)\Big)P_{d_H}(d\gamma),
 	\end{multline*}
 	where $d_{\rho_1}(M)=\max_{x,y\in M} \rho_1(x,y)$.
Let $\gamma\in[0,d_{\rho_1}(M)]$ be fixed,
 	\begin{multline*}
 		\int_{M^2} \mathbb{P}(x\in B_1(\partial A_p(x_1,x_2),\gamma))P_X(dx_1) P_X(dx_2)=\\
 		P_X\otimes P_X\{(x_1,x_2)\in M^2:x\in B_1(\partial A_p(x_1,x_2),\gamma) \}=  \\
 		\int_M P_X\{x_1\in M:x\in B_1(\partial A_p(x_1,x_2),\gamma)\}P_X(dx_2).
 	\end{multline*}
 	Let  $f_0=\max_{x\in M}f(x)$, and $\nu$ be the volume measure on $M$ inherited from $g_1$, then, for fixed $\gamma$ and $x_2$,
  $P_X\{x_1\in M:x\in B_1(\partial A_p(x_1,x_2),\gamma)\}= P_X\Big[ B_1(\partial B_p(x_1,\rho_p(x,x_2)),\gamma)\Big]\leq  f_0 \nu(B_1(\partial B_p(x_1,\rho_p(x,x_2)),\gamma)).$ 
 	From  the dominated convergence theorem it follows easily that for fixed $\gamma$ and $x$, the function $x_2\to \nu(B_1(\partial B_p(x_1,\rho_p(x,x_2)),\gamma))$ is a continuous function of $x_2$. Its maximum is reached at some $x_2^0=x_2^0(\gamma)$.
Let us bound  $\nu(B_1(\partial B_p(x_1,\rho_p(x,x^0_2)),\gamma))$.
By condition I, $\partial B_p(x_1,\rho_p(x,x^0_2))$  has positive reach w.r.t. $\rho_p$. Then by the Corollary on page 57  in \cite{bangert1982}, it has positive reach $R>0$ w.r.t. $\rho_1$.
So that, there exists a positive constant $\tau$ such that
$\nu(B_1(\partial B_p(x_1,\rho_p(x,x^0_2)),\gamma))=\tau\gamma \quad  \text{ for all } \gamma\in [0,R].$ 
 	This result was initially proved in \cite{federer59} for subsets of $\mathbb{R}^d$ and generalized to manifolds in \cite{klein81}.
 	
 	Then $I$ is bounded from above by
 	$$\beta_n\delta^{-1}f_0\tau \int_0^{R} \gamma P_{d_H}(d\gamma)+ 
 		\beta_n\delta^{-1}f_0 \int_R^{d_{\rho_1}(M)}  \nu(B_1(\partial B_p(x_1,\rho_p(x,x_2^0)),\gamma))P_{d_H}(d\gamma)=J_1+J_2.$$
	To bound $J_2$ multiply and divide the integral by $\gamma$, and bound $1/\gamma\leq 1/R$. So that,  $J_2$ is bounded from above by
 	\begin{equation}\label{eq0}
 		\beta_n\delta^{-1}f_0 \nu(M)R^{-1} \int_R^{d_{\rho_1}(M)}   \gamma P_{d_H}(d\gamma)\leq \beta_n\delta^{-1}f_0 \nu(M)R^{-1} \int_0^{d_{\rho_1}(M)} P(d_H>\gamma)d\gamma.
 	\end{equation}
 	
 	Now $J_1$ is bounded from above by 
 	\begin{equation}\label{eq2}
 		\beta_n\delta^{-1}f_0\tau \int_0^{d_{\rho_1}(M)} \gamma P_{d_H}(d\gamma)=f_0\beta_n 2L_0\delta^{-1}\tau \int_0^{d_{\rho_1}(M)} P (d_H>\gamma)d\gamma.
 	\end{equation}
 	 Corollary \ref{c1} implies that \eqref{eq0} and \eqref{eq2} are bounded from above by $C\delta^{-1}\beta_n \exp(-\theta_0 n^{1/(d+2p)}).$
 	We may conclude that \eqref{rate} holds using  $n^{1/(d+2p)}/(\log(\beta_n)\log(n))\to \infty$ and the Borel--Cantelli Lemma.
 \QEDB


\end{document}